\documentclass[12pt,twoside]{amsart}
\usepackage[T1]{fontenc}
\usepackage{amsmath,amsthm,amsfonts,amssymb}
\usepackage[bookmarksnumbered,plainpages=false,colorlinks=true,linkcolor=blue,citecolor=blue,urlcolor=blue]{hyperref}
\numberwithin{equation}{section}
\newtheorem{theorem}{Theorem}[section]
\newtheorem{proposition}[theorem]{Proposition}
\newtheorem{lemma}[theorem]{Lemma}
\newtheorem{corollary}[theorem]{Corollary}
\theoremstyle{definition}
\newtheorem{definition}[theorem]{Definition}
\newtheorem{remark}[theorem]{Remark}
\newtheorem{example}[theorem]{Example}
\newcommand{\Res}{\operatorname{Wres}}
\newcommand{\Rtheta}{\operatorname{Res}_{\theta}}
\newcommand{\tr}{\operatorname{tr}}
\newcommand{\Tr}{\operatorname{Tr}}
\newcommand{\Ric}{\operatorname{Ric}}
\newcommand{\vol}{\operatorname{dvol}_{g}}
\newcommand{\id}{\mathbf{1}}
\newcommand{\nc}{\mathcal{A}_{\theta}}

\newcommand{\Om}{\Omega}
\newcommand{\dd}{\partial}
\newcommand{\sphere}{\nu_{n-1}}
\newcommand{\Cn}{C_n}
\newcommand{\Kn}{K_n}
\newcommand{\owedge}{\mathbin{\wedge}}
\title{Clifford-weighted residues, partial spectral curvature, and torsion}
\author{Yong Wang}
\address{School of Mathematics and Statistics, Northeast Normal University,
Changchun, Jilin 130024, China}
\email{wangy581@nenu.edu.cn}
\subjclass[2020]{58J42; 58B34; 53C27; 46L87}
\keywords{Noncommutative residue, Dirac operator, partial spectral curvature,
Riemann curvature, noncommutative torus, torsion}
\date{}
\begin{document}
\begin{abstract}
We study Clifford-weighted noncommutative residues associated with Dirac operators and connections. On closed spin manifolds a double commutator of the square of the Dirac operator, paired with a Clifford commutator, recovers the localized Riemann curvature tensor. We prove presentation independence for its extension to arbitrary classical one-forms, recompute four-Clifford functionals, and give explicit Weyl, Schouten and Einstein curvature combinations. For spectral triples with an admissible residue calculus we define a partial spectral curvature functional on two arbitrary represented one-forms and two exact represented one-forms. Its well-definedness is distinguished from module tensoriality: explicit flat noncommutative-torus computations exhibit both a failure of general presentation descent and nonzero partial pairings. A connection-based curvature pairing provides a separate construction, including noncentral pure Clifford-degree-two examples. We also derive ordered first- and second-derivative torsion residues on free torus modules and prove an obstruction to ordinary alternation for unrestricted noncentral inputs. All constants, coefficient orders and hypotheses are kept explicit.
\end{abstract}
\maketitle
\pagestyle{myheadings}
\markboth{Yong Wang}{Clifford-weighted residues and partial spectral curvature}
\begingroup\small\tableofcontents\endgroup
\medskip

\section{Introduction}
Spectral residues encode geometric information through insertions into inverse powers of elliptic operators. The metric and Einstein functionals of D\k{a}browski, Sitarz and Zalecki \cite{DSZmetric}, their spectral torsion functional \cite{DSZtorsion}, and the general Laplace-type insertion calculation of Bochniak, D\k{a}browski, Sitarz and Zalecki \cite{BDSZ} form the analytic background. Their general symbol formulas are used here with stated normalizations. Our purpose is to separate full-curvature pairings, Ricci-level four-Clifford data, and the extra terms arising from noncommutative coefficients.

For a closed spin manifold, consider
\begin{equation}
 \mathcal U_D(a_0;f_1,f_2,f_3,f_4)
 =\Res\bigl(a_0[[D,f_1],[D,f_2]]
                [[D^2,f_3],[D^2,f_4]]D^{-n}\bigr).
\end{equation}
The square commutator separates into a first-order covariant derivative, spin curvature and scalar multiplication. The compatible first-order residue vanishes, and the fiber trace of the Clifford commutator eliminates the scalar term. The remaining curvature term is
$-4K_n\int a_0R(\operatorname{grad}f_3,\operatorname{grad}f_4,\operatorname{grad}f_1,\operatorname{grad}f_2)\vol$.
We give the complete product-rule and Clifford calculations and prove independence of classical one-form presentations. Arbitrary test functions recover the local density, not just one integrated component.

The covariant insertions motivating the construction are
\begin{align}
 \mathcal F(X,Y,Z,W)&=\Res\bigl(c(X)c(Y)\nabla_Z^S\nabla_W^SD^{-n}\bigr),\label{eq:Fintro}\\
 \mathcal H_T(X,Y,Z)&=\Res\bigl(c(X)c(Y)\nabla_Z^{T,S}D_T^{-n}\bigr).\label{eq:Hintro}
\end{align}
The first contains scalar, Ricci and Riemann curvature contributions. Antisymmetrization extracts full curvature. The second detects the difference between the inserted torsion connection and the Bochner connection of $D_T^2$. It is proportional to the established spectral torsion functional, rather than a separate torsion invariant.

For general spectral triples we restrict the last two inputs to exact represented forms and define a partial spectral curvature pairing directly through $[D^2,f]$. This avoids lifting a chosen sum representing a general one-form. We prove that the general presentation lift need not descend, even on the flat noncommutative torus.  These results identify well-defined operator functionals while keeping their tensor properties separate.

Connection-based noncommutative Riemann geometry requires additional module, metric and compatibility structure; see Mesland and Rennie \cite{MRLC,MRcurv}. We do not identify an arbitrary spin-degree-two bundle curvature with a Levi-Civita Riemann tensor. Instead we give tested fixed-frame curvature coefficients, a noncentral example, and sufficient hypotheses for agreement with the classical construction. The Hodge--Dirac setting \cite{DSZHodge} provides a further commutative-algebra descent example, whereas a direct sum with a finite matrix triple supplies a noncommutative algebra whose finite sector is invisible to residues.

We recompute the four-Clifford functionals appearing in \cite{LW}, fixing specifically the supplied version arXiv:2503.01084v1 in comparisons. Their corrected expressions depend only on scalar and Ricci curvature. Combining their antisymmetric parts with $\mathcal U_D$ gives explicit Weyl and related curvature formulas. This comparison concerns displayed coefficients, not a claim that a single numerical residue reconstructs a metric.

The noncommutative-torus calculations are performed on free modules with flat scalar principal symbol, using the calculus of \cite{HLP1,HLP2,Ponge}. Every product of coefficients is retained in order. The resulting connection and torsion formulas are exact in that setting; noncentral curved principal symbols, projective-module calculi and boundary residues require separate analysis. The claims are confined to the formulas and hypotheses proved here.

\section{Conventions and the SIGMA insertion formula}
\subsection{Spin and residue conventions}
Let $(M,g)$ be a closed Riemannian spin manifold of dimension $n=2m\geq2$. Use the full complex spinor bundle $S$, of rank $r=2^m$, and
\begin{equation}
 c(u)c(v)+c(v)c(u)=-2g(u,v)\id,\qquad
 D=\sum_a c(e_a)\nabla^S_{e_a}.
\end{equation}
For a differential form $\alpha$,
\begin{equation}
 c(\alpha)=\frac1{k!}\alpha_{a_1\dots a_k}c(e_{a_1})\cdots c(e_{a_k}),\quad \alpha\in\Omega^k(M).
\end{equation}
Our curvature conventions are
\begin{align}
 R(Z,W)&=[\nabla^{LC}_Z,\nabla^{LC}_W]-\nabla^{LC}_{[Z,W]},\\
 R(Z,W,X,Y)&=g(R(Z,W)X,Y),\\
 \Ric(Z,W)&=\sum_a g(R(e_a,Z)W,e_a),\qquad s=\tr_g\Ric.
\end{align}
Thus the round unit sphere has positive scalar curvature, whereas $R(e_1,e_2,e_1,e_2)=-1$ for an orthonormal two-plane on that sphere. The spin curvature is
\begin{equation}
 \Om^S(Z,W)=\frac14\sum_{a,b}g(R(Z,W)e_a,e_b)c(e_a)c(e_b).
 \label{eq:spinOmega}
\end{equation}
We use the normalization
\begin{equation}
 \Res(P)=(2\pi)^{-n}\int_M\int_{|\xi|=1}\tr\sigma_{-n}(P)(x,\xi)\,dS(\xi)\,\vol.
 \label{eq:resnormal}
\end{equation}
Set
\begin{equation}
 \sphere=\operatorname{Vol}(S^{n-1})=\frac{2\pi^{n/2}}{\Gamma(n/2)},\quad
 \Cn=(2\pi)^{-n}\sphere,\quad \Kn=r\Cn.
 \label{eq:constants}
\end{equation}
Our residue is $(2\pi)^{-n}$ times the unnormalized cosphere integral used in \cite{BDSZ}. In that paper the generators $\gamma^a$ have positive Clifford square; the conversion is $c(e_a)=i\gamma^a$. All inverse powers below are parametrices modulo smoothing operators. Consequently a nontrivial kernel does not affect any formula.

\subsection{The insertion formula}
In normal coordinates at a point, suppose a Laplace-type operator has symbol
\begin{align}
 p_2&=\left(\delta_{ab}+\tfrac13 \mathcal R_{acbd}x^cx^d\right)\xi_a\xi_b+o(|x|^2),\\
 p_1&=i(P_{ab}x^b+S_a)\xi_a+o(|x|),\qquad p_0=Q+o(1),
\end{align}
where $\mathcal R_{acbd}=-R(e_a,e_c,e_b,e_d)$ in our four-slot convention. This distinction keeps the principal metric expansion in the convention of \cite{BDSZ} while retaining \eqref{eq:spinOmega}. Suppose also that at the chosen point
\begin{equation}
 \sigma(O)=F^{ab}\xi_a\xi_b+iG^a\xi_a+H,\qquad F^{ab}=F^{ba}.
\end{equation}
Proposition 2.1 of \cite{BDSZ}, in normalization \eqref{eq:resnormal}, reads
\begin{align}
 \Res(OL^{-m})=\frac{\Cn}{24}\int_M\tr\bigl[&24H+12G^aS_a\nonumber\\
 &+F^{aa}(-12Q+6P_{bb}-2s\id-3S_bS_b)\nonumber\\
 &+2F^{ab}(-6P_{ab}+2\Ric_{ab}\id-3S_aS_b)\bigr]\,\vol.
 \label{eq:SIGMA}
\end{align}
Only the symmetric part of $P_{ab}$ enters its contraction with $F^{ab}$. We will use a radial bundle frame to make the connection contributions transparent.

\subsection{Sphere moments}
\begin{lemma}[Sphere moments]\label{lem:annSphere}
Writing $\sphere=\operatorname{Vol}(S^{n-1})$,
\begin{align}
 \int\xi_i\xi_j\,dS&=\frac{\sphere}{n}\delta_{ij},\label{eq:annMoment2}\\
 \int\xi_i\xi_j\xi_k\xi_l\,dS
 &=\frac{\sphere}{n(n+2)}
 (\delta_{ij}\delta_{kl}+\delta_{ik}\delta_{jl}+\delta_{il}\delta_{jk}).
 \label{eq:annMoment4}
\end{align}
\end{lemma}

\subsection{Bochner data and the covariant insertion}

\begin{lemma}[Radial-frame data]\label{lem:radial}
Let $L=\nabla^{E*}\nabla^E+V$ act on a Hermitian bundle $E$. In normal coordinates and a radial parallel frame,
\begin{equation}
 \nabla^E_a=\dd_a+A_a(x),\qquad
 A_a(x)=\tfrac12x^b\Om^E_{ba}+O(|x|^2).
\end{equation}
The data in \eqref{eq:SIGMA} satisfy
\begin{equation}
 S_a=0,\quad Q=V,\quad
 P_{ab}=\tfrac23\Ric_{ab}\id+\Om^E_{ab},\quad
 P_{aa}=\tfrac23s\id.
 \label{eq:radialdata}
\end{equation}
\end{lemma}
\begin{proof}
The radial-gauge identity follows by differentiating radial parallel transport; at the origin $A_a=0$ and $\dd_bA_a=\Om^E_{ba}/2$. The first-order coordinate coefficient of the scalar connection Laplacian is $\frac23\Ric_{ab}x^b\dd_a$. The bundle connection adds $-2A_a\dd_a=-x^b\Om^E_{ba}\dd_a=x^b\Om^E_{ab}\dd_a$. Its zero-order connection contribution is $-\sum_a\dd_aA_a=0$ at the origin. The zeroth-order coefficient is therefore $V$; the antisymmetry of $\Om^E$ gives the last identity.
\end{proof}

\subsubsection{Radial parallel frames and the complete insertion calculation}\label{sec:radialannotations}
We give a detailed derivation of every coefficient in Lemma \ref{lem:radial}, including the choice of frame and the derivatives which remain nonzero at its center. Throughout this subsubsection $p$ is fixed, $x(p)=0$, and all coefficients without an explicit argument in a Taylor expansion are evaluated at $p$. The letter $A_i$ denotes a connection matrix; the separate endomorphism multiplying an insertion on the left will be denoted $B$.

\paragraph{Construction of the radial parallel frame.}
Choose a sufficiently small normal coordinate ball about $p$. For $x$ in this ball the curve
\begin{equation}
 \gamma_x(t)=\exp_p(t x),\qquad 0\leq t\leq1,
\end{equation}
has coordinate representation $t x$. Select a basis $v_1,\ldots,v_N$ of $E_p$ and define a local frame $s_\alpha(x)$ by parallel transport of $v_\alpha$ along $\gamma_x$ using $\nabla^E$. Smooth dependence of the parallel-transport ordinary differential equation on its endpoint gives a smooth frame. Uniqueness and composition of parallel transport show that, for fixed $x$, $s_\alpha(t x)$ is parallel along this same radial curve. Thus
\begin{equation}
 0=\nabla^E_{\dot\gamma_x(t)}s_\alpha(t x)
   =x^i\nabla^E_{\partial_i}s_\alpha(t x).
\end{equation}
Define the matrix $A_i$ by
\begin{equation}
 \nabla^E_{\partial_i}s_\alpha
 =\sum_\beta(A_i)^\beta{}_{\alpha}s_\beta.
\end{equation}
Since the frame vectors are independent, $x^iA_i(t x)=0$, in particular
\begin{equation}
 \boxed{x^iA_i(x)=0.}
 \label{eq:radialgaugeexplicit}
\end{equation}
For a column of component functions $\psi$, the same definition and the connection product rule give
\begin{equation}
 \nabla_i^E\psi=(\partial_i+A_i)\psi.
\end{equation}
The normal coordinate choice concerns the metric on the base; radial parallel transport concerns the bundle connection. They can therefore be imposed simultaneously, without assuming that $\nabla^E$ is the Levi-Civita connection.

\paragraph{Why the connection matrices also vanish at the center.}
Substitution of $x=0$ into \eqref{eq:radialgaugeexplicit} alone gives only $0=0$. Instead differentiate the identity, which holds on the entire coordinate ball:
\begin{equation}
 0=\partial_k(x^iA_i)=A_k+x^i\partial_kA_i.
\end{equation}
At the center this implies
\begin{equation}
 \boxed{A_k(0)=0\quad\text{for every }k.}
 \label{eq:radialcenterzero}
\end{equation}
Equivalently, on the $k$th coordinate axis $x=t e_k$ with $t\ne0$, the gauge identity implies $A_k(t e_k)=0$; continuity gives the value at the origin. This does not imply that $A_k(x)$ is zero off that axis, nor that its derivatives vanish at the origin.

\paragraph{Recovering the first derivative from curvature.}
In coordinate directions the curvature matrix is
\begin{equation}
 \Omega^E_{ki}=\partial_kA_i-\partial_iA_k+[A_k,A_i].
\end{equation}
Differentiate the radial condition with respect to $x^i$ to obtain $x^k\partial_iA_k=-A_i$. Also $x^k[A_k,A_i]=[x^kA_k,A_i]=0$. Consequently
\begin{equation}
 x^k\Omega^E_{ki}=x^k\partial_kA_i+A_i.
 \label{eq:radialcurvaturecontraction}
\end{equation}
Along a fixed ray, the chain rule gives
\begin{align}
 \frac{d}{dt}\bigl(tA_i(t x)\bigr)
 &=A_i(t x)+t x^k\partial_kA_i(t x)\nonumber\\
 &=t x^k\Omega^E_{ki}(t x).
\end{align}
The term $tA_i(t x)$ tends to zero as $t\to0$, because the connection matrices are smooth. Integrating from zero to one proves the exact identity
\begin{equation}
 \boxed{A_i(x)=\int_0^1t x^k\Omega^E_{ki}(t x)\,dt.}
 \label{eq:radialIntegral}
\end{equation}
Since $\Omega^E_{ki}(t x)=\Omega^E_{ki}(0)+O(|x|)$ uniformly for $0\leq t\leq1$,
\begin{equation}
 A_i(x)=\frac12x^k\Omega^E_{ki}(0)+O(|x|^2).
\end{equation}
In particular,
\begin{equation}
 \boxed{\partial_kA_i(0)=\frac12\Omega^E_{ki}(0),
 \qquad\sum_i\partial_iA_i(0)=0.}
 \label{eq:radialderivativeexplicit}
\end{equation}
The second equality uses $\Omega^E_{ii}=0$. One can check the factor $1/2$ without integration: twice differentiating $x^kA_k=0$ at zero yields $\partial_kA_i+\partial_iA_k=0$. Together with $A_i(0)=0$ in the curvature formula, this gives $\Omega^E_{ki}(0)=2\partial_kA_i(0)$.

Nonzero curvature is fully compatible with \eqref{eq:radialcenterzero}. For example, on a trivial Hermitian line bundle in a plane, set
\begin{equation}
 A_1=-\frac{i\lambda}{2}x^2,\qquad
 A_2=\frac{i\lambda}{2}x^1,\qquad\lambda\in\mathbb R\setminus\{0\}.
\end{equation}
Then $x^1A_1+x^2A_2=0$ and $A_1(0)=A_2(0)=0$, but $\Omega^E_{12}=i\lambda$. Curvature is retained in the first derivatives, not in the values of these matrices at the center.

\paragraph{The scalar first-order coefficient in normal coordinates.}
Let $v(x)=\sqrt{\det g(x)}$. For the positive scalar Laplacian,
\begin{align}
 \Delta f
 &=-v^{-1}\partial_i(vg^{ij}\partial_jf)\nonumber\\
 &=-g^{ij}\partial_i\partial_jf
  -\bigl(\partial_i g^{ij}+g^{ij}\partial_i\log v\bigr)\partial_jf.
 \label{eq:radialscalarlaplacian}
\end{align}
To track its signs, use the auxiliary curvature components
$\mathcal R_{ikjl}=-R(e_i,e_k,e_j,e_l)$ already specified before \eqref{eq:SIGMA}. The normal metric expansions in this convention are
\begin{align}
 g_{ij}(x)&=\delta_{ij}-\frac13\mathcal R_{ikjl}x^kx^l+O(|x|^3),\nonumber\\
 g^{ij}(x)&=\delta_{ij}+\frac13\mathcal R_{ikjl}x^kx^l+O(|x|^3).
\end{align}
The inverse expansion follows from multiplying the two matrices: products of their quadratic corrections have order four. Since $\sum_i\mathcal R_{ikil}=\Ric_{kl}$, the determinant identity $\det(I+H)=1+\operatorname{tr}H+O(H^2)$ yields
\begin{equation}
 v(x)=1-\frac16\Ric_{kl}x^kx^l+O(|x|^3),\qquad
 \log v(x)=-\frac16\Ric_{kl}x^kx^l+O(|x|^3).
\end{equation}
Differentiation, together with symmetry of Ricci, gives
\begin{equation}
 \partial_j\log v=-\frac13\Ric_{jk}x^k+O(|x|^2).
\end{equation}
For the inverse metric divergence,
\begin{align}
 \sum_i\partial_i g^{ij}
 &=\frac13\sum_{i,k,l}\mathcal R_{ikjl}
       (\delta_i^k x^l+x^k\delta_i^l)+O(|x|^2)\nonumber\\
 &=\frac13\sum_{i,k}\mathcal R_{ikji}x^k+O(|x|^2)
 =-\frac13\Ric_{jk}x^k+O(|x|^2).
\end{align}
Here $\mathcal R_{iijl}=0$ removes the first contraction, and $\sum_i\mathcal R_{ikji}=-\Ric_{kj}$ gives the second. Thus both terms inside the parentheses in \eqref{eq:radialscalarlaplacian} contribute $-\Ric_{jk}x^k/3$. Including the outer minus sign, the scalar Laplacian has first-order coefficient
\begin{equation}
 \frac23\Ric_{jk}x^k\partial_j+O(|x|^2)\partial_j.
\end{equation}

\paragraph{The bundle first-order and zero-order coefficients.}
The connection Laplace operator has the local expression
\begin{equation}
 L=-g^{ij}\bigl((\partial_i+A_i)(\partial_j+A_j)
              -\Gamma_{ij}^k(\partial_k+A_k)\bigr)+V.
\end{equation}
Its first-order bundle contribution, to linear order in $x$, is
\begin{equation}
 -2A_j\partial_j
 =-x^k\Omega^E_{kj}\partial_j+O(|x|^2)\partial_j
 =x^k\Omega^E_{jk}\partial_j+O(|x|^2)\partial_j.
\end{equation}
Combining with the scalar coefficient gives the symbol data
\begin{equation}
 S_j=0,\qquad P_{jk}=\frac23\Ric_{jk}\id+\Omega^E_{jk},
 \qquad\sum_jP_{jj}=\frac23s\id.
\end{equation}
At zero the zeroth-order bundle expression is
$-\sum_i\partial_iA_i-\sum_iA_i^2+\sum_{i,k}\Gamma_{ii}^kA_k$.
It vanishes by \eqref{eq:radialcenterzero}, \eqref{eq:radialderivativeexplicit} and the normal-coordinate identity $\Gamma_{ij}^k(0)=0$. Hence $Q=V(0)$. This proves every entry of \eqref{eq:radialdata}; $Q$ is the Bochner potential at this point, not an additional spin trace or a multiple of scalar curvature.

\paragraph{Full expansion of the two-derivative insertion.}
Write $Z=Z^i\partial_i$ and $W=W^j\partial_j$. For a section represented by a column $\psi$,
\begin{equation}
 \nabla_Z^E\nabla_W^E\psi
 =Z^i(\partial_i+A_i)(W^j\partial_j\psi+W^jA_j\psi).
\end{equation}
There are two product rules:
\begin{align}
 \partial_i(W^j\partial_j\psi)
 &=(\partial_iW^j)\partial_j\psi+W^j\partial_i\partial_j\psi,\nonumber\\
 \partial_i(W^jA_j\psi)
 &=(\partial_iW^j)A_j\psi+W^j(\partial_iA_j)\psi
                          +W^jA_j\partial_i\psi.
\end{align}
The matrix multiplication by $A_i$ adds two further terms. Therefore the exact expansion is
\begin{equation}
 \begin{split}
 \nabla_Z^E\nabla_W^E\psi=Z^i\bigl(&W^j\partial_i\partial_j\psi
 +(\partial_iW^j)\partial_j\psi+W^jA_j\partial_i\psi
 +A_iW^j\partial_j\psi\\
 &+(\partial_iW^j)A_j\psi+W^j(\partial_iA_j)\psi
 +A_iW^jA_j\psi\bigr).
 \end{split}
 \label{eq:radialfullinsertion}
\end{equation}
Thus all seven terms are accounted for: two from the first product rule, three from the second, and two from left multiplication by $A_i$. There is no derivative of $Z^i$, since it is the outer left coefficient. Derivatives of $W^j$ do occur. There is no missing Christoffel correction: this is the composition $\nabla_Z^E(\nabla_W^E\psi)$, not the tensorial Hessian
$\nabla_Z^E\nabla_W^E\psi-\nabla^E_{\nabla_Z^{LC}W}\psi$.

At the center all terms with an undifferentiated $A_i$ vanish, but $\partial_iA_j(0)=\Omega^E_{ij}(0)/2$ need not vanish. Consequently
\begin{equation}
 \begin{split}
 (\nabla_Z^E\nabla_W^E\psi)(0)
 ={}&Z^iW^j\partial_i\partial_j\psi
   +Z^i(\partial_iW^j)\partial_j\psi\\
 &+\frac12Z^iW^j\Omega^E_{ij}\psi,
 \end{split}
\end{equation}
where every factor on the right is evaluated at zero. Left multiplication by a smooth endomorphism $B$ and the convention $\partial_i\mapsto i\xi_i$ now give
\begin{equation}
 \boxed{
 F^{ij}=-B Z^{(i}W^{j)},\qquad
 G^j=B Z^i\partial_iW^j,\qquad
 H=\frac12 B Z^iW^j\Omega^E_{ij}.}
\end{equation}
Here $Z^{(i}W^{j)}=(Z^iW^j+Z^jW^i)/2$; symmetrization is allowed because $\xi_i\xi_j$ is symmetric. The minus sign in $F$ comes from $\partial_i\partial_j\mapsto-\xi_i\xi_j$. The curvature term in $H$ is the surviving derivative of the connection matrix. It must not be dropped merely because the matrices themselves vanish at the center.

Finally, inserting these coefficients into \eqref{eq:SIGMA} kills $G$ because $S=0$, not because $\partial_iW^j$ is zero. It also kills the contraction of the symmetric $F^{ij}$ with the antisymmetric curvature part of $P_{ij}$. Curvature nevertheless remains through $H$. The left multiplier $B$ is not differentiated by this composition rule. These observations explain precisely how the residue in Proposition \ref{prop:two} becomes tensorial even though the differential-operator composition itself is not.

\begin{remark}[Antisymmetric connection curvature]
The curvature part of $P_{ij}$ has zero contraction with the symmetric $F^{ij}$. This does not remove curvature from the answer: $H$ independently contains $\Om(Z,W)/2$. Omitting that zero-order insertion coefficient would incorrectly leave only scalar and Ricci terms in the two-derivative functional.
\end{remark}

\begin{proposition}[Two-derivative insertion]\label{prop:two}
For any smooth endomorphism $A$ and vector fields $Z,W$,
\begin{align}
 \Res(A\nabla^E_Z\nabla^E_WL^{-m})
 =\Cn\int_M\tr\biggl[A\biggl(&\tfrac16\Ric(Z,W)\id
 -\tfrac1{12}s\,g(Z,W)\id\nonumber\\
 &+\tfrac12g(Z,W)V+\tfrac12\Om^E(Z,W)\biggr)\biggr]\,\vol.
 \label{eq:generalmanifold}
\end{align}
Moreover $\Res(A\nabla^E_UL^{-m})=0$ for every $U$.
\end{proposition}
\begin{proof}
At the origin of the radial frame the coefficients of $O=A\nabla_Z^E\nabla_W^E$ are
\begin{align}
 F^{ab}&=-A Z^{(a}W^{b)},\qquad F^{aa}=-A g(Z,W),\\
 G^a&=A Z^b\dd_bW^a,\qquad
 H=\tfrac12 A Z^aW^b\Om^E_{ab}.
\end{align}
The term $G^aS_a$ is zero. By Lemma \ref{lem:radial},
\begin{align}
 -12Q+6P_{bb}-2s\id&=-12V+2s\id,\\
 -6P_{ab}+2\Ric_{ab}\id&=-2\Ric_{ab}\id-6\Om^E_{ab}.
\end{align}
The last antisymmetric term contracts to zero with $F^{ab}$. Substitution in \eqref{eq:SIGMA} gives, successively, the contributions
\begin{equation}
 \tfrac12 A\Om^E(Z,W),\quad
 A g(Z,W)(\tfrac12V-\tfrac1{12}s\id),\quad
 \tfrac16 A\Ric(Z,W),
\end{equation}
which prove \eqref{eq:generalmanifold}. For $O=A\nabla_U^E$, one has $F=H=0$ at the origin and $S=0$, so \eqref{eq:SIGMA} gives zero. No derivative of $A$ occurs because it is a left multiplication insertion.
\end{proof}

\begin{remark}[Tensoriality after taking the residue]
The composition $\nabla_Z^E\nabla_W^E$ is not $C^\infty(M)$-linear in $W$. The residue is unchanged on replacing it by
\begin{equation}
 \nabla_Z^E\nabla_W^E-\nabla^E_{\nabla^{LC}_ZW},
\end{equation}
since the correction has zero residue by Proposition \ref{prop:two}. Thus \eqref{eq:generalmanifold} is a tensor density in $Z,W$. This property is specific to an insertion using the Bochner connection of $L$.
\end{remark}

\section{Covariant insertions and Riemann curvature}\label{sec:covariant}
\begin{lemma}[Clifford traces]\label{lem:trace}
For the full spin representation,
\begin{align}
 \tr_S(c(X)c(Y))&=-r g(X,Y),\\
 \tr_S(c(X)c(Y)c(e_a)c(e_b))
 &=r\bigl(g(X,Y)\delta_{ab}-X_aY_b+X_bY_a\bigr),\label{eq:fourtrace}\\
 \tr_S(c(X)c(Y)\Om^S(Z,W))&=-\tfrac r2R(Z,W,X,Y),\label{eq:Rtrace}\\
 \tr_S(c(X)c(Y)c(\iota_ZT))&=-rT(Z,X,Y).
 \label{eq:Ttrace}
\end{align}
\end{lemma}
\begin{proof}
The Clifford relation and cyclicity give the two-generator identity and \eqref{eq:fourtrace}. Contracting \eqref{eq:fourtrace} with the antisymmetric coefficient in \eqref{eq:spinOmega} leaves $\frac r4(-R(Z,W,X,Y)+R(Z,W,Y,X))$, proving \eqref{eq:Rtrace}. For \eqref{eq:Ttrace}, use $c(\iota_ZT)=\frac12T(Z,e_a,e_b)c(e_a)c(e_b)$ and the same contraction. The factor $1/2$ in the two-form Clifford multiplication is essential.
\end{proof}

\subsection{Clifford trace identities explained}\label{sec:traceannotations}
The traces above are taken in the full complex spin representation, not a half-spin representation. Here is a direct derivation that also fixes the signs of the torsion contractions.

For an ordered index set $I=\{i_1<\cdots<i_p\}$ write $c_I=c_{i_1}\cdots c_{i_p}$. The Clifford relations imply
\begin{equation}
 c_I^2=(-1)^{p(p+1)/2}\id.
 \label{eq:bladeSquare}
\end{equation}
There are $p(p-1)/2$ interchanges to pair equal generators, followed by $p$ factors $c_i^2=-1$. If $I$ is nonempty, then $\tr(c_I)=0$: for even $p$ choose a generator inside $I$, and for odd $p$ choose one outside $I$. Conjugation by that invertible generator changes the sign of $c_I$, while preserving trace. The outside generator exists for odd $p$ because $n$ is even. In particular the top-degree blade also has zero full spinor trace. Distinct blades have zero trace pairing, since their product is a nonempty blade up to sign. Therefore
\begin{equation}
 \tr(c(\alpha)c(\beta))
 =r(-1)^{p(p+1)/2}\langle\alpha,\beta\rangle,
 \quad \alpha,\beta\in\Omega^p(M).
 \label{eq:gradeTrace}
\end{equation}
For different pure degrees the trace pairing is zero. In increasing-index notation our convention reads
\begin{equation}
 c(\alpha)=\sum_{i_1<\cdots<i_p}\alpha_{i_1\dots i_p}c_{i_1}\cdots c_{i_p}.
\end{equation}
This is equivalent to the factorial convention stated above: each increasing-index term occurs $p!$ times in the unrestricted sum, with the signs of the alternating coefficient and of the Clifford permutation canceling.

For completeness cyclicity gives the following recursion for any even number of vector factors:
\begin{align}
 \tr(c(v_1)\cdots c(v_{2k}))
 =\sum_{j=2}^{2k}(-1)^{j-1}g(v_1,v_j)
 \tr(c(v_2)\cdots\widehat{c(v_j)}\cdots c(v_{2k})).
 \label{eq:traceRecursion}
\end{align}
Move $c(v_1)$ through all other factors. The final cyclic trace is the negative of the original one, because $2k-1$ is odd. The anticommutators contribute $-2g$, and division by two gives \eqref{eq:traceRecursion}. For $k=1$ it is $-rg(v_1,v_2)$; for $k=2$ it gives the three metric pairings in \eqref{eq:fourtrace}.

Spin curvature is a Clifford two-form. Formula \eqref{eq:gradeTrace} for $p=2$ makes its trace pairing negative. More explicitly, insert \eqref{eq:spinOmega} into \eqref{eq:fourtrace}; the $g(X,Y)\delta_{ab}$ term vanishes because $R(Z,W,e_a,e_a)=0$. The two remaining terms give
\[
 \frac r4\bigl[-R(Z,W,X,Y)+R(Z,W,Y,X)\bigr]
 =-\frac r2R(Z,W,X,Y).
\]
Replacing $c(X)c(Y)$ by its commutator doubles this value. This is the factor used in Theorem \ref{thm:DiracR}.

Likewise $c(\iota_ZT)$ is a two-form, so
\begin{equation}
 \tr(c(X)c(Y)c(\iota_ZT))=-rT(Z,X,Y).
\end{equation}
By contrast the pure degree-three trace pairing has positive sign, since $(-1)^{3\cdot4/2}=1$. The product of three vector factors decomposes as
\begin{equation}
 c(X)c(Y)c(Z)=c(X^\flat\wedge Y^\flat\wedge Z^\flat)
 -g(X,Y)c(Z)+g(X,Z)c(Y)-g(Y,Z)c(X).
\end{equation}
The degree-one terms pair to zero with $c(T)$, and hence
\begin{equation}
 \tr(c(X)c(Y)c(Z)c(T))=rT(X,Y,Z).
 \label{eq:annSixTrace}
\end{equation}
The different signs in these two torsion traces are thus required by Clifford degree; they do not signal a change of torsion convention.

\begin{remark}[What changes with noncommuting coefficients]
These vector trace identities use scalar component functions. Over a noncommutative coefficient algebra one may still use constant Clifford matrix traces, but the coefficients remain in order. For instance $\tr_S(c(X)c(Y))=-r\sum_i x^iy^i$, whereas reversing the arguments gives $-r\sum_i y^ix^i$. The two expressions need not be equal as algebra elements. The torus formulas below retain precisely this distinction.
\end{remark}

\subsection{Deriving the Lichnerowicz potential and evaluating the residue}
At a synchronous point Clifford compatibility gives
\begin{equation}
 D^2=-\sum_i(\nabla_i^S)^2+\frac12\sum_{i,j}c_ic_j\Om^S_{ij}.
\end{equation}
The symmetric part of $c_ic_j$ produces the connection Laplacian; the antisymmetric part produces the commutator of covariant derivatives. The contracted spin identity \eqref{eq:contractedspin}, proved directly from Bianchi and Clifford algebra below, gives
\begin{equation}
 \sum_{i,j}c_ic_j\Om^S_{ij}
 =-\frac12\sum_i c_i c(\Ric^\sharp e_i)=\frac{s}{2}\id.
\end{equation}
Here symmetry of the Ricci matrix kills the degree-two part of $\sum_i c_ic(\Ric^\sharp e_i)$, leaving $-s\id$. Restoring the coordinate correction in the connection Laplacian yields the global identity $D^2=\nabla^{S*}\nabla^S+s\id/4$. The scalar contribution to the insertion density is consequently $-s g(Z,W)/12+s g(Z,W)/8=s g(Z,W)/24$, which accounts for the coefficient in the next theorem.

\begin{theorem}[Evaluation of the first functional]\label{thm:classicalF}
With the conventions above,
\begin{align}
 \mathcal F(X,Y,Z,W)=-\Kn\int_M\biggl[&g(X,Y)
 \left(\tfrac16\Ric(Z,W)+\tfrac1{24}s\,g(Z,W)\right)\nonumber\\
 &+\tfrac14R(Z,W,X,Y)\biggr]\,\vol.
 \label{eq:Ffinal}
\end{align}
\end{theorem}
\begin{proof}
The Lichnerowicz identity gives $D^2=\nabla^{S*}\nabla^S+\frac14s\id$. Thus $V=s\id/4$ in \eqref{eq:generalmanifold}, and
\begin{equation}
 -\frac{s}{12}g(Z,W)\id+\frac12g(Z,W)V
 =\frac{s}{24}g(Z,W)\id.
\end{equation}
Take $A=c(X)c(Y)$ and use Lemma \ref{lem:trace}. The scalar endomorphism terms yield the first line of \eqref{eq:Ffinal}, and the spin-curvature term yields its second line.
\end{proof}

\begin{corollary}[Curvature extraction]\label{cor:curvature}
The antisymmetric part in either pair of arguments is
\begin{equation}
 \frac{\mathcal F(X,Y,Z,W)-\mathcal F(Y,X,Z,W)}2
 =-\frac{\Kn}{4}\int_M R(Z,W,X,Y)\,\vol.
 \label{eq:antisym}
\end{equation}
Equivalently,
\begin{equation}
 \Res\bigl(c(X)c(Y)\Om^S(Z,W)D^{-n}\bigr)
 =-\frac{\Kn}{2}\int_M R(Z,W,X,Y)\,\vol.
\end{equation}
\end{corollary}
\begin{proof}
The metric terms are symmetric in $X,Y$ and the curvature term is antisymmetric. Alternatively, antisymmetrize the derivative composition, use its commutator identity, and discard $\nabla^S_{[Z,W]}$ by Proposition \ref{prop:two}. The curvature endomorphism times $D^{-n}$ contributes only its principal symbol.
\end{proof}
\noindent These statements concern integrated pairings. Inserting an arbitrary scalar test function $f$ to the left multiplies the displayed density by $f$. Therefore the collection of localized pairings determines the pointwise curvature tensor, not merely an integral of its components.

\subsection{Contractions and Einstein correction}
Removing the Clifford insertion gives
\begin{equation}
 \Res(\nabla^S_Z\nabla^S_WD^{-n})
 =\Kn\int_M\left(\tfrac16\Ric(Z,W)+\tfrac1{24}s\,g(Z,W)\right)\vol.
 \label{eq:Ricplus}
\end{equation}
For $n>2$, the usual scalar-curvature residue is
\begin{equation}
 \Res(fD^{-(n-2)})=-\frac{\Kn(n-2)}{24}\int_M f s\,\vol.
 \label{eq:scalar}
\end{equation}
This follows from \cite[Corollary 2.2]{BDSZ}, or from the heat coefficient $a_2(D^2)=-r s/12$ and the residue--heat relation. Consequently
\begin{align}
 \Res\left[\left(\nabla^S_Z\nabla^S_W+
 \frac{3g(Z,W)}{n-2}D^2\right)D^{-n}\right]
 =\frac{\Kn}{6}\int_M\left(\Ric-\tfrac12s g\right)(Z,W)\,\vol.
 \label{eq:Einstein}
\end{align}
The restriction $n>2$ is necessary for this correction. Formula \eqref{eq:Ffinal}, on the other hand, is valid also in dimension two.
\begin{remark}[Deriving the Einstein correction coefficient]
For a scalar $a$, replacing the correction coefficient in \eqref{eq:Einstein} by $a/(n-2)$ changes the scalar part of the density to
\begin{equation}
 \Kn\left(\frac{s}{24}-\frac{as}{24}\right)g(Z,W)
 =\frac{\Kn}{24}(1-a)s\,g(Z,W).
\end{equation}
To obtain $\Kn(\Ric-sg/2)/6$ this coefficient must be $-\Kn s/12$, so $1-a=-2$ and $a=3$. The choice $a=1$ cancels scalar curvature entirely and yields only $\Kn\Ric/6$. The factor three corrects the coefficient omitted in the preceding draft; it is fixed here independently of any literature normalization.
\end{remark}

\subsection{Checks of signs and normalizations}
On a flat torus the right-hand side of \eqref{eq:Ffinal} vanishes. For a unit sphere, $\Ric=(n-1)g$ and $s=n(n-1)$. At an orthonormal two-plane,
\begin{equation}
 \mathcal F(e_1,e_2,e_1,e_2)\text{ has density }\Kn/4,
\end{equation}
in our curvature convention. Also summing \eqref{eq:Ricplus} over an orthonormal frame gives $\Kn(n+4)s/24$. Directly,
\begin{equation}
 \sum_a\left(\nabla^S_{e_a}\nabla^S_{e_a}
 -\nabla^S_{\nabla^{LC}_{e_a}e_a}\right)
 =-D^2+s/4,
\end{equation}
and \eqref{eq:scalar} together with $\Res(fD^{-n})=\Kn\int_M f\vol$ gives exactly the same coefficient. These checks also detect a possible incorrect scalar-curvature sign.

\section{Riemann curvature from Dirac double commutators}\label{sec:double}
The curvature insertion in Corollary \ref{cor:curvature} uses a connection explicitly. We now eliminate that input: the operator in the next theorem involves only $D$, multiplication by functions, and commutators. Throughout this section $\Delta=-\operatorname{div}\operatorname{grad}$ is the positive scalar Laplacian.

\begin{lemma}[Commutators with the square]\label{lem:squarecomm}
For $f\in C^\infty(M)$,
\begin{equation}
 [D,f]=c(df),\qquad [D^2,f]=-2\nabla^S_{\operatorname{grad}f}+(\Delta f)\id.
 \label{eq:squarecomm}
\end{equation}
If $Z=\operatorname{grad}f$ and $W=\operatorname{grad}h$, then
\begin{align}
 [[D^2,f],[D^2,h]]
 &=4\nabla^S_{[Z,W]}+4\Om^S(Z,W)+q(f,h)\id,\label{eq:doubleexpand}\\
 q(f,h)&=-2Z(\Delta h)+2W(\Delta f).
\end{align}
\end{lemma}
\begin{proof}
The first identity follows from the Leibniz rule for $D$. At a point with a synchronous orthonormal frame, $\nabla^{S*}\nabla^S=-\sum_i\nabla_i^S\nabla_i^S$ up to the usual connection correction. Applying the product rule twice yields
\[
 [\nabla^{S*}\nabla^S,f]\psi=(\Delta f)\psi-2\nabla^S_{\operatorname{grad}f}\psi.
\]
The potential $s/4$ in $D^2$ commutes with $f$, proving \eqref{eq:squarecomm}. Write the two operators as $-2\nabla_Z^S+b$ and $-2\nabla_W^S+d$, with $b=\Delta f$, $d=\Delta h$. Their commutator equals
\[
 4[\nabla_Z^S,\nabla_W^S]-2Z(d)\id+2W(b)\id.
\]
Substitute $[\nabla_Z^S,\nabla_W^S]=\nabla_{[Z,W]}^S+\Om^S(Z,W)$. In particular this commutator has order at most one, because the first-order principal coefficients are scalar.
\end{proof}

\begin{theorem}[A Dirac-only Riemann functional]\label{thm:DiracR}
For $a_0,f_1,f_2,f_3,f_4\in C^\infty(M)$ define
\begin{align}
 \mathcal U_D(a_0;f_1,f_2,f_3,f_4)
 =\Res\Big(&a_0[[D,f_1],[D,f_2]]\nonumber\\
 &\cdot[[D^2,f_3],[D^2,f_4]]D^{-n}\Big).
 \label{eq:Udefinition}
\end{align}
Then
\begin{align}
 \mathcal U_D(a_0;f_1,f_2,f_3,f_4)
 =-4\Kn\int_M a_0 R(&\operatorname{grad}f_3,\operatorname{grad}f_4,\nonumber\\
 &\operatorname{grad}f_1,\operatorname{grad}f_2)\,\vol.
 \label{eq:Uformula}
\end{align}
\end{theorem}
\begin{proof}
Set
\begin{equation}
 \begin{aligned}
 X&=\operatorname{grad}f_1,&Y&=\operatorname{grad}f_2,\\
 Z&=\operatorname{grad}f_3,&W&=\operatorname{grad}f_4,
 \end{aligned}
 \qquad B=[c(X),c(Y)].
\end{equation}
Lemma \ref{lem:squarecomm} splits \eqref{eq:Udefinition} into three residues. Proposition \ref{prop:two} makes the residue of the first-order insertion vanish, including derivatives of the variable bracket field. Since $\tr B=0$, the scalar term contributes
\[
 \Res(a_0BqD^{-n})=\Cn\int_M a_0q\tr B\,\vol=0.
\]
Here only the principal symbol $|\xi|^{-n}\id$ contributes: this insertion has order zero, and lower symbols have degree below $-n$. For the remaining term Lemma \ref{lem:trace} gives
\begin{align*}
 \tr\bigl(B\Om^S(Z,W)\bigr)
 &=-\frac r2R(Z,W,X,Y)+\frac r2R(Z,W,Y,X)\\
 &=-rR(Z,W,X,Y).
\end{align*}
Consequently its residue is $4\Cn\int_M a_0\tr(B\Om^S(Z,W))\vol$, proving \eqref{eq:Uformula} with the factor $-4\Kn$.
\end{proof}

\subsection{Arbitrary one-forms and independence of presentation}
Every smooth one-form on a compact manifold admits a finite presentation $\alpha=\sum_j a_j\,db_j$. Indeed, choose a finite coordinate cover and a subordinate partition of unity with supports inside the charts. Extend each coordinate function after multiplying by a cutoff that is one on the relevant support; extending the supported coefficients by zero gives a global finite presentation.
For such a presentation set
\begin{equation}
 L_\alpha=\sum_j a_j[D^2,b_j],\qquad
 c(\alpha)=\sum_j a_j[D,b_j].
\end{equation}
Clifford multiplication is intrinsic, whereas $L_\alpha$ can depend on the presentation. More precisely,
\begin{equation}
 L_\alpha=-2\nabla^S_{\alpha^\sharp}+b_\alpha\id,
 \qquad b_\alpha=\sum_j a_j\Delta b_j.
 \label{eq:Lalpha}
\end{equation}

\begin{corollary}[The four-one-form version]\label{cor:oneforms}
For arbitrary presentations of $\alpha,\beta,\gamma,\delta$,
\begin{align}
 \Res\bigl(a_0[c(\alpha),c(\beta)][L_\gamma,L_\delta]D^{-n}\bigr)
 =-4\Kn\int_M a_0R(\gamma^\sharp,\delta^\sharp,
 \alpha^\sharp,\beta^\sharp)\,\vol.
 \label{eq:Uoneforms}
\end{align}
The left side is independent of the chosen presentations.
\end{corollary}
\begin{proof}
Formula \eqref{eq:Lalpha} follows from \eqref{eq:squarecomm}; the coefficients $a_j$ multiply on the left, so they are not differentiated in this formula. Expand $[L_\gamma,L_\delta]$ exactly as in Lemma \ref{lem:squarecomm}. It equals $4\nabla^S_{[\gamma^\sharp,\delta^\sharp]}+4\Om^S(\gamma^\sharp,\delta^\sharp)$ plus a scalar endomorphism. The preceding proof applies verbatim. Alternatively, a change of presentation changes $L_\gamma$ and $L_\delta$ by scalar multiplications. Their commutator changes by a scalar multiplication, whose product with $[c(\alpha),c(\beta)]$ has zero spinor trace.
\end{proof}

\subsection{Tensoriality, localization and recovery}
The numerical functional is an integrated tensor pairing, not by itself a pointwise tensor. Define its normalized local density by
\begin{equation}
 \mathfrak R_D(\alpha,\beta,\gamma,\delta)
 =-\frac1{4\Kn}\operatorname{dens}_g
 \Res\bigl([c(\alpha),c(\beta)][L_\gamma,L_\delta]D^{-n}\bigr),
\end{equation}
where $\operatorname{dens}_g$ means the scalar coefficient relative to $\vol$. The proof above is pointwise in a radial frame and shows
\begin{equation}
 \mathfrak R_D(\alpha,\beta,\gamma,\delta)
 =R(\gamma^\sharp,\delta^\sharp,\alpha^\sharp,\beta^\sharp).
 \label{eq:localR}
\end{equation}
Thus this density is $C^\infty(M)$-multilinear, antisymmetric in each pair, symmetric under pair exchange, and satisfies the first Bianchi identity. These are conclusions of the calculation; individual operators $L_\alpha$ need not be module-linear in their presentation.

For recovery, fix a chart and extend its coordinate functions smoothly using cutoffs. On any smaller region where the cutoffs are one, their differentials form a coframe. Allowing $a_0$ to range over test functions supported in that region determines each smooth component of \eqref{eq:localR}: if its integral against every such $a_0$ is zero, the component vanishes. Pair symmetry identifies this slot order with $R(\alpha^\sharp,\beta^\sharp,\gamma^\sharp,\delta^\sharp)$. A single unlocalized value of $\mathcal U_D$ does not contain this recovery information.

\begin{remark}[Dimension and orientation]
The double-commutator argument also works in dimension two. The four-dimensional lower bound is needed below for the usual Weyl and Schouten decompositions. No chirality operator is inserted here: the full spinor trace is essential to the displayed constants and avoids orientation-dependent characteristic terms.
\end{remark}

\subsection{Product rules and the exact curvature isolation}\label{sec:classicalproducts}
The second identity in Lemma \ref{lem:squarecomm} can be checked without discarding any connection correction. In a local orthonormal frame let $V_a=\nabla^{LC}_{e_a}e_a$. Then
\begin{align}
 \nabla^{S*}\nabla^S&=-\sum_a(\nabla^S_{e_a}\nabla^S_{e_a}-\nabla^S_{V_a}),\nonumber\\
 \Delta f&=-\sum_a(e_ae_af-V_af).
\end{align}
For every spinor $\psi$ the two product rules are
\begin{align}
 \nabla^S_{e_a}\nabla^S_{e_a}(f\psi)
 &=(e_ae_af)\psi+2(e_af)\nabla^S_{e_a}\psi
          +f\nabla^S_{e_a}\nabla^S_{e_a}\psi,\nonumber\\
 \nabla^S_{V_a}(f\psi)&=(V_af)\psi+f\nabla^S_{V_a}\psi.
\end{align}
Subtracting $f\nabla^{S*}\nabla^S\psi$ gives $(\Delta f)\psi-2\nabla^S_{\operatorname{grad}f}\psi$. The Lichnerowicz potential $s/4$ commutes with $f$. This proves the lemma globally, rather than only in a special frame.

There is a useful Clifford verification. The algebraic identity
$[D^2,f]=D[D,f]+[D,f]D$ and Clifford compatibility give
\begin{equation}
 [D^2,f]=c_a c(\nabla^{LC}_{e_a}\operatorname{grad}f)
                   -2\nabla^S_{\operatorname{grad}f}.
\end{equation}
The Hessian $H_{ab}=g(\nabla_{e_a}^{LC}\operatorname{grad}f,e_b)$ is symmetric, so
$H_{ab}c_ac_b=H_{ab}(c_ac_b+c_bc_a)/2=-\sum_aH_{aa}\id=(\Delta f)\id$.
This also fixes the sign for the positive scalar Laplacian.

For the double commutator write $P=-2\nabla_Z+b$, $Q=-2\nabla_W+d$. Their complete products on $\psi$ are
\begin{align}
 PQ\psi={}&4\nabla_Z\nabla_W\psi-2Z(d)\psi
             -2d\nabla_Z\psi-2b\nabla_W\psi+bd\psi,\nonumber\\
 QP\psi={}&4\nabla_W\nabla_Z\psi-2W(b)\psi
             -2b\nabla_W\psi-2d\nabla_Z\psi+db\psi.
\end{align}
The mixed first-order terms cancel and $bd=db$, leaving
$4[\nabla_Z,\nabla_W]-2Z(d)+2W(b)$. No assertion that gradient fields commute is used. For example $Z=x\partial_x$ and $W=y\partial_x+x\partial_y$ are gradients of $x^2/2$ and $xy$, but $[Z,W]=-y\partial_x+x\partial_y$.

The proof of Theorem \ref{thm:DiracR} involves three separate vanishing mechanisms. The bracket-field insertion has zero residue because it uses the Bochner connection of $D^2$. The scalar multiplication term is killed by the fiber trace of the Clifford commutator. The remaining zero-order curvature insertion is evaluated using only the leading inverse symbol, which supplies $\Cn$. These mechanisms yield $4(-r)\Cn=-4\Kn$ and remain valid with a variable left test coefficient $a_0$.

\subsection{Global one-form presentations and scalar lift errors}\label{sec:globalpresentations}
Here is a precise globalization of the finite presentation used in Corollary \ref{cor:oneforms}. Choose finitely many coordinate charts $U_\lambda$ and a partition of unity with $\operatorname{supp}\rho_\lambda\Subset U_\lambda$. In a chart write $\alpha=\alpha_{\lambda k}dx_\lambda^k$. The initially global identity is
\begin{equation}
 \alpha=\sum_\lambda\rho_\lambda\alpha.
\end{equation}
Each local one-form $\rho_\lambda\alpha_{\lambda k}dx_\lambda^k$ vanishes near the chart boundary, so has a smooth zero extension $\eta_{\lambda k}$; globally $\rho_\lambda\alpha=\sum_k\eta_{\lambda k}$. Coordinate differentials themselves are not being extended without a cutoff.

Choose $\chi_\lambda\in C_c^\infty(U_\lambda)$ equal to one on an open neighborhood of $\operatorname{supp}\rho_\lambda$. Extend
$a_{\lambda k}=\rho_\lambda\alpha_{\lambda k}$ and
$b_{\lambda k}=\chi_\lambda x_\lambda^k$ by zero. The exact calculation in the chart is
\begin{equation}
 a_{\lambda k}db_{\lambda k}
 =\rho_\lambda\alpha_{\lambda k}\chi_\lambda dx_\lambda^k
  +\rho_\lambda\alpha_{\lambda k}x_\lambda^k d\chi_\lambda
 =\eta_{\lambda k}.
\end{equation}
Indeed $\rho_\lambda\chi_\lambda=\rho_\lambda$ and $\rho_\lambda d\chi_\lambda=0$. Where $\rho_\lambda$ can be nonzero, $\chi_\lambda$ is constant on a neighborhood, not merely equal to one at a single point. Outside the chart both zero extensions vanish. Consequently $\alpha=\sum_{\lambda,k}a_{\lambda k}db_{\lambda k}$ as a global finite sum.

For any such presentation,
\begin{equation}
 \sum_j a_j[D,b_j]=c(\alpha),\qquad
 L_\alpha=-2\nabla^S_{\alpha^\sharp}+u_\alpha\id,
 \qquad u_\alpha=\sum_j a_j\Delta b_j.
\end{equation}
Thus the Clifford action is intrinsic, whereas the scalar $u_\alpha$ can change. For arbitrary $\gamma,\delta$, put $Z=\gamma^\sharp$, $W=\delta^\sharp$, $u=u_\gamma$, $v=u_\delta$. The same product calculation above gives
\begin{equation}
 [L_\gamma,L_\delta]
 =4\nabla^S_{[Z,W]}+4\Omega^S(Z,W)
           +(-2Z(v)+2W(u))\id.
\end{equation}
Hence the proof of Theorem \ref{thm:DiracR} applies verbatim and proves Corollary \ref{cor:oneforms} for nonexact forms. More directly, if presentations change by $L_\gamma'=L_\gamma+e$, $L_\delta'=L_\delta+f$, then
\begin{equation}
 [L_\gamma',L_\delta']-[L_\gamma,L_\delta]
 =(-2Z(f)+2W(e))\id.
\end{equation}
Its product with $[c(\alpha),c(\beta)]$ has zero fiber trace and hence zero residue. This is a residue-level independence theorem; it does not assert that $L_\alpha$ itself is intrinsic. Nor can one prove it simply by pulling all presentation coefficients outside commutators, since a first-order operator differentiates such coefficients.

\section{Partial spectral curvature and descent of one-form presentations}\label{sec:partial}
The preceding manifold construction has two distinct extensions. A connection supplies a curvature endomorphism directly. Alternatively, the exact commutators with $D^2$ define an intrinsic operation on exact represented one-forms. The latter requires no choice of connection, but its tensor properties must be distinguished from its well-definedness. In this section $\mathfrak f$ denotes a specified residue trace; its normalization is not changed during a calculation.

\subsection{Analytic assumptions and exact one-forms}
Let $(\mathcal A,\mathcal H,D)$ be a unital spectral triple with a smooth pseudodifferential calculus and a residue trace defined on all products appearing below. We assume that $\mathcal A$ and $[D,\mathcal A]$ preserve a common smooth core for $D$, that the compositions with $D^2$ are defined there, and that their products with $|D|^{-n}$ belong to the residue calculus. These are assumptions on the calculus, not consequences of compact resolvent alone. In the classical and flat-torus examples they follow from the ordinary or toroidal calculus. Replace $|D|$ on its kernel by an invertible smoothing modification whenever necessary.

The represented one-forms and their exact subspace are
\begin{equation}
 \Omega_D^1(\mathcal A)=\left\{\sum_j a_j[D,b_j]\right\},\qquad
 \mathcal E_D^1=\{[D,f]:f\in\mathcal A\}.
 \label{eq:partialspaces}
\end{equation}
Both are complex vector spaces; the first is an $\mathcal A$-bimodule, whereas the second is generally not an $\mathcal A$-submodule. Define, on the common core,
\begin{equation}
 \Lambda_D([D,f])=[D^2,f]=\{D,[D,f]\}.
 \label{eq:exactlift}
\end{equation}
This is well defined on the represented exact form. Indeed $[D,f]=[D,g]$ implies
$[D^2,f-g]=D[D,f-g]+[D,f-g]D=0$.

\begin{definition}[Partial spectral Riemann curvature functional]\label{def:partial}
For $a_0\in\mathcal A$, $\alpha,\beta\in\Omega_D^1(\mathcal A)$ and
$\eta,\zeta\in\mathcal E_D^1$, set
\begin{equation}
 \mathcal P_D(a_0;\alpha,\beta;\eta,\zeta)
 =-\frac14\mathfrak f\left(a_0[\alpha,\beta]
 [\Lambda_D(\eta),\Lambda_D(\zeta)]|D|^{-n}\right).
 \label{eq:partialdefinition}
\end{equation}
The term ``partial'' specifies the restriction of the last two slots to exact represented one-forms. We also call this a partial spectral curvature pairing. The word tensor will be used only when additional module and curvature identities have been established.
\end{definition}

\begin{proposition}[Intrinsic slots and classical normalization]\label{prop:partialclassical}
The functional \eqref{eq:partialdefinition} is independent of all presentations of its four inputs. It is complex multilinear in the four form slots and antisymmetric in each displayed pair. For the canonical spin triple and $\mathfrak f=\Res$,
\begin{equation}
 \mathcal P_D(a_0;c(\alpha),c(\beta);[D,f],[D,h])
 =\Kn\int_M a_0R(\operatorname{grad}f,\operatorname{grad}h,
                         \alpha^\sharp,\beta^\sharp)\vol.
 \label{eq:partialclassical}
\end{equation}
\end{proposition}
\begin{proof}
The first two forms enter as operators, so equal represented forms give equal commutators. The last two enter through the well-defined map \eqref{eq:exactlift}. Linearity of $\Lambda_D$, of the commutator in each slot, and of $\mathfrak f$ proves complex multilinearity. Reversing either pair changes the sign of its commutator. In the spin case use Lemma \ref{lem:squarecomm}, with $Z=\operatorname{grad}f$, $W=\operatorname{grad}h$, to split the last commutator into
$4\nabla^S_{[Z,W]}+4\Omega^S(Z,W)+q\id$.
The first term has zero residue by Proposition \ref{prop:two}. The scalar term is killed by the pointwise identity $\tr_S[c(\alpha),c(\beta)]=0$. Lemma \ref{lem:trace} gives
$\tr_S([c(\alpha),c(\beta)]\Omega^S(Z,W))=-rR(Z,W,\alpha^\sharp,\beta^\sharp)$.
The factor $-1/4$ in the definition cancels the factor $4$ in the double commutator and its trace sign. The principal inverse symbol and sphere integral give $\Cn$, proving \eqref{eq:partialclassical}.
\end{proof}

No universal spin rank exists for an arbitrary spectral triple; thus $\Kn$ is not put into the general definition. On the spin triple, $\mathcal P_D/\Kn$ is the normalized integrated Riemann pairing. Pair exchange and Bianchi identities follow there from the recovered curvature density. They do not follow merely from the two algebraic commutators in Definition \ref{def:partial}. In particular, multiplication of an exact form by an arbitrary $a\in\mathcal A$ may leave $\mathcal E_D^1$, so unrestricted module linearity in its last two slots is not even a property on that domain.

\subsection{Why the presentation lift of a general form is different}
For a chosen expression $\omega=\sum a_j[D,b_j]$, write
$L_\omega=\sum a_j[D^2,b_j]$. Moving $D$ past each left coefficient, without changing any coefficient order, gives
\begin{align}
 L_\omega
 &=\sum_j a_jD[D,b_j]+\sum_j a_j[D,b_j]D\nonumber\\
 &=D\omega+\omega D-\sum_j[D,a_j][D,b_j].
 \label{eq:liftobstruction}
\end{align}
Consequently a null presentation $\sum a_j[D,b_j]=0$ can have the nonzero lift
\begin{equation}
 E=\sum_j a_j[D^2,b_j]=-\sum_j[D,a_j][D,b_j].
 \label{eq:nullpresentationerror}
\end{equation}
Let $\mathcal J_D$ be the complex space of all such lift errors. For fixed represented forms $\gamma,\delta$ and chosen lifts, changing them by $E,F\in\mathcal J_D$ changes their commutator by
\begin{equation}
 [L_\gamma+E,L_\delta+F]-[L_\gamma,L_\delta]
 =[E,L_\delta]+[L_\gamma,F]+[E,F].
 \label{eq:descenterror}
\end{equation}
Thus the presentation-based functional descends to represented one-forms exactly when
\begin{equation}
 \mathfrak f\left(a_0[\alpha,\beta]
 ([E,L_\delta]+[L_\gamma,F]+[E,F])|D|^{-n}\right)=0
 \label{eq:descentcondition}
\end{equation}
for all allowed forms, presentations and errors. Residue traciality alone does not establish this: formally, for admissible products,
$\mathfrak f(A[E,L]P)=\mathfrak f(E[L,PA])$, which need not vanish.

\begin{proposition}[A commutative Dirac-type descent criterion]\label{prop:clifforddescent}
For the scalar algebra $C^\infty(M)$ represented on a Hermitian Clifford module and a self-adjoint Dirac-type operator $D$, condition \eqref{eq:descentcondition} holds for the Wodzicki residue. This includes Hodge--Dirac and spin$^c$ triples as well as twisted spin Dirac triples with scalar function algebra.
\end{proposition}
\begin{proof}
Write $D^2=\nabla^{B*}\nabla^B+V$. The product-rule proof of Lemma \ref{lem:squarecomm} applies to this Laplace-type operator and gives
$[D^2,f]=-2\nabla^B_{\operatorname{grad}f}+(\Delta f)\id$.
The Clifford symbol identifies a represented one-form with its ordinary one-form, so
$L_\omega=-2\nabla^B_{\omega^\sharp}+u_\omega\id$.
A null presentation therefore has lift $e\id$. Hence $E=e\id$, $F=f\id$, $[E,F]=0$, and the other two commutators in \eqref{eq:descenterror} are scalar multiplication, namely $2\delta^\sharp(e)\id$ and $-2\gamma^\sharp(f)\id$. A scalar multiple of $[\alpha,\beta]$ has zero fiber trace. The zero-order residue formula proves \eqref{eq:descentcondition}. This argument proves descent; it does not identify the resulting Clifford-module curvature trace with a spin Riemann tensor without an additional calculation.
\end{proof}

The scalar-algebra hypothesis is essential in this proof. Replacing it by matrix-valued functions introduces noncentral coefficients in the same infinite-dimensional residue. Spectral closedness in \cite[Definition 5.5 and Theorem 5.7]{DSZmetric} concerns a different vanishing and a module-balancing property of their Einstein functional. It cannot be substituted for \eqref{eq:descentcondition}. The Hodge--Dirac spectral functionals are studied in \cite{DSZHodge}.

\begin{example}[A noncommutative direct sum satisfying descent]\label{ex:directsum}
Let $(C^\infty(M),\mathcal H_M,D_M)$ satisfy Proposition \ref{prop:clifforddescent}. Set
\begin{equation}
 \mathcal A=C^\infty(M)\oplus M_N(\mathbb C),\quad
 \mathcal H=\mathcal H_M\oplus\mathbb C^N,\quad
 D=D_M\oplus D_F,\qquad D_F=\operatorname{diag}(1,\ldots,N).
\end{equation}
This is a spectral triple with block-diagonal representation. The finite commutators need not vanish: $[D_F,E_{ij}]=(i-j)E_{ij}$. Both represented forms and their presentation lifts split into the two blocks. For a block operator $P=P_M\oplus P_F$,
\begin{equation}
 \Tr(P|D|^{-z})=\Tr(P_M|D_M|^{-z})
                 +\sum_{j=1}^N(P_F)_{jj}j^{-z}.
 \label{eq:directsumzeta}
\end{equation}
The finite sum is entire, since $j^{-z}=\exp(-z\log j)$, independently of whether $P_F$ commutes with $D_F$. Thus the zeta residue equals that of the classical block. Every error in \eqref{eq:descentcondition} also splits, and the classical block vanishes by Proposition \ref{prop:clifforddescent}; the finite block has zero residue. Hence descent holds for all represented forms of the direct sum. The residue itself can be nonzero on other inputs, but it detects none of the finite noncommutative sector. This example is a direct sum, not the interacting product algebra $C^\infty(M)\otimes M_N(\mathbb C)$.
\end{example}

\subsection{A flat noncommutative torus where presentation descent fails}\label{sec:descentfailure}
Use the torus conventions introduced in Section \ref{sec:toruscalculus}. Let two generators satisfy $UV=qVU$ with $|q|=1$, $q\ne1$, and
$\partial_1U=iU$, $\partial_1V=0$, $\partial_2U=0$, $\partial_2V=iV$.
Set $D_0=c_i\partial_i$, $r=2^{n/2}$ and use $\Rtheta$ as the residue. Then
\begin{equation}
 U^*[D_0,U]+U[D_0,U^*]=0,
 \label{eq:torusnullform}
\end{equation}
whereas
\begin{align}
 [D_0^2,U]&=-2iU\partial_1+U,&
 [D_0^2,U^*]&=2iU^*\partial_1+U^*,\nonumber\\
 U^*[D_0^2,U]+U[D_0^2,U^*]&=2\id.
 \label{eq:torusnulllift}
\end{align}
Multiplying the null presentation by $t/2$ on the left shows that $\mathcal J_{D_0}$ contains every multiplication operator $t\in\nc$.

Take $\alpha=Uc_1$, $\beta=Vc_1$, $\gamma=0$, and $\delta=c_1=-iU^*[D_0,U]$. With this presentation,
$L_\delta=-2\partial_1-i$. Choose initially $L_\gamma=0$ and then change its zero-form presentation so that $L_\gamma'=t$. The other error is zero. Since
$[\alpha,\beta]=-[U,V]\id$ and $[t,L_\delta]=2\partial_1t$, the change of the unnormalized presentation functional is
\begin{equation}
 \Delta=-2r\sphere\tau([U,V]\partial_1t).
\end{equation}
Put $t=(UV)^*$. Then $\partial_1t=-it$ and
$[U,V]=(1-q^{-1})UV$, giving
\begin{equation}
 \Delta=2ir\sphere(1-q^{-1})\ne0.
 \label{eq:descentcounterexample}
\end{equation}
Thus \eqref{eq:descentcondition} fails on this standard flat triple. On the restricted domain where the first two represented forms have central coefficients, their spin commutator has zero spin trace. A null lift here is a Clifford scalar: if $\sum a_j\partial_kb_j=0$, applying $\partial_i$ and subtracting the $i,k$ equations gives
$\sum(\partial_i a_j)(\partial_kb_j)-(\partial_k a_j)(\partial_i b_j)=0$; hence the degree-two part of \eqref{eq:nullpresentationerror} vanishes. Error commutators then have only spin-scalar coefficients, and pair trivially with that restricted spin-traceless factor in the flat residue. This restricted cancellation does not contradict \eqref{eq:descentcounterexample}.

\subsection{A partial pairing can be nonzero on the flat triple}
Definition \ref{def:partial} has no presentation problem in the preceding example. Nevertheless, flatness alone does not force it to vanish. Take
$\alpha=Uc_1$, $\beta=Vc_1$, $\eta=[D_0,U]$, $\zeta=[D_0,V]$, and put $s=[U,V]$. The first two inputs are represented forms because $c_1=-iU^*[D_0,U]$. Direct multiplication gives
\begin{align}
 [\alpha,\beta]&=-s\id,\nonumber\\
 [[D_0^2,U],[D_0^2,V]]
 &=s(-4\partial_1\partial_2-2i\partial_1-2i\partial_2+1).
 \label{eq:partialflatoperators}
\end{align}
The flat inverse symbol is $|\xi|^{-n}\id$ modulo smoothing and has no algebra derivatives. The order-two and order-one terms therefore have no homogeneous component of degree $-n$ after multiplication by $|D_0|^{-n}$. Only the zero-order coefficient contributes. It follows that
\begin{equation}
 \mathcal P_{D_0}(a_0;\alpha,\beta;\eta,\zeta)
   =\frac{r\sphere}{4}\tau(a_0s^2).
 \label{eq:partialflatvalue}
\end{equation}
Taking $a_0=(s^2)^*$ yields the positive value
$r\sphere|1-q^{-1}|^4/4$. Indeed $s=(1-q^{-1})UV$ is a nonzero scalar multiple of a unitary, so $s^2$ has squared norm $|1-q^{-1}|^4$. This is a localized partial spectral pairing on a flat triple, not the curvature of its flat frame connection. It demonstrates that exactness of the last slots removes presentation ambiguity but not the additional noncommutative terms.

\section{Four Clifford insertions and their curvature content}\label{sec:PQ}
We examine the two functionals of Li and Wang \cite[Theorem 1.1]{LW}, specifically the supplied version arXiv:2503.01084v1. Their cosphere residue omits $(2\pi)^{-n}$. Every comparison below first converts their formulas to \eqref{eq:resnormal}; in particular their constant $2^m\nu_{n-1}$ becomes $\Kn$.
For vector fields $X,Y,Z,W$ put
\begin{align}
 P_D(X,Y,Z,W)&=\Res(c(X)c(Y)c(Z)c(W)D^{-(n-2)}),\\
 Q_D(X,Y,Z,W)&=\Res(c(X)c(Y)D c(Z)c(W)D D^{-n}),
 \label{eq:PQdef}
\end{align}
and write
\begin{align}
 J(X,Y,Z,W)&=g(X,Z)g(Y,W)-g(X,W)g(Y,Z),\\
 (g\owedge H)(X,Y,Z,W)&=g(X,Z)H(Y,W)+g(Y,W)H(X,Z)\nonumber\\
 &\quad-g(X,W)H(Y,Z)-g(Y,Z)H(X,W)
 \label{eq:KNproduct}
\end{align}
for a symmetric two-tensor $H$. Notice $g\owedge g=2J$.

\subsection{Clifford contractions needed for the recomputation}
\begin{lemma}\label{lem:PQclifford}
Let $C=c(Z)c(W)$ and let $\Ric^\sharp$ denote the Ricci endomorphism. In an orthonormal frame, with $c_i=c(e_i)$,
\begin{align}
 \sum_i c_i Cc_i&=(4-n)C+4g(Z,W)\id,\label{eq:sandwich}\\
 \sum_{i,j}\Ric_{ij}c_i Cc_j
 &=-sC-2c(W)c(\Ric^\sharp Z)+2c(Z)c(\Ric^\sharp W),\label{eq:Ricsandwich}\\
 \sum_{i,j}c_i Cc_j\Om^S_{ij}
 &=\frac{s}{2}C+c(W)c(\Ric^\sharp Z)-c(Z)c(\Ric^\sharp W).
 \label{eq:Omegasandwich}
\end{align}
\end{lemma}
\begin{proof}
Move $c_i$ through the two factors of $C$ using the Clifford relation. More explicitly,
\begin{equation}
 c_i C=Cc_i+2W_i c(Z)-2Z_i c(W).
 \label{eq:moveCi}
\end{equation}
Summing after right multiplication by $c_i$ gives $-nC+2c(Z)c(W)-2c(W)c(Z)$, which is \eqref{eq:sandwich}. Contracting \eqref{eq:moveCi} with $\Ric_{ij}c_j$ instead gives \eqref{eq:Ricsandwich}, because symmetry of Ricci implies $\sum\Ric_{ij}c_ic_j=-s\id$.

The contracted spin-curvature identity in the present convention is
\begin{equation}
 \sum_j c_j\Om^S_{ij}=-\frac12c(\Ric^\sharp e_i).
 \label{eq:contractedspin}
\end{equation}
For completeness insert \eqref{eq:spinOmega}. The degree-three part of $\sum_{j,a,b}R(e_i,e_j,e_a,e_b)c_jc_ac_b/4$ vanishes by pair symmetry and the first Bianchi identity. The repeated-index terms leave
\[
 \frac12\sum_{j,b}R(e_i,e_j,e_b,e_j)c_b
 =-\frac12\sum_b\Ric_{ib}c_b,
\]
which proves \eqref{eq:contractedspin}. Consequently $\sum_{i,j}c_ic_j\Om^S_{ij}=s\id/2$. Apply \eqref{eq:moveCi} again and then \eqref{eq:contractedspin}: the term $C\sum c_ic_j\Om^S_{ij}$ is $sC/2$, and the two remaining terms are $-c(Z)c(\Ric^\sharp W)+c(W)c(\Ric^\sharp Z)$. This proves \eqref{eq:Omegasandwich}.
\end{proof}

\begin{remark}
The sandwich trace in the specified version of \cite[Eq. (3.23)]{LW} differs from \eqref{eq:sandwich}. In particular, when all four Clifford vectors are the same unit vector, the actual summed trace is $-nr$, rather than $+nr$. The recomputation below uses the Clifford identity directly.
\end{remark}

\subsection{The contracted spin-curvature identity with all indices}\label{sec:contractedspindetail}
The last identity in Lemma \ref{lem:PQclifford} uses
\begin{equation}
 \sum_jc_j\Omega^S_{ij}=-\frac12c(\Ric^\sharp e_i).
 \label{eq:contractedspinexpanded}
\end{equation}
To prove it, write $R_{ijab}=R(e_i,e_j,e_a,e_b)$ and insert
$\Omega^S_{ij}=R_{ijab}c_ac_b/4$. The elementary Clifford decomposition is
\begin{equation}
 c_jc_ac_b=c(e^j\wedge e^a\wedge e^b)
             -\delta_{ja}c_b+\delta_{jb}c_a-\delta_{ab}c_j.
\end{equation}
For distinct indices this is the definition of the three-form action. When two indices coincide it follows from $c_j^2=-1$ and anticommutation; if all coincide, the right side is $-c_j+c_j-c_j=-c_j=c_j^3$.

The three-form sum $\mathcal T_i=\sum_{j,a,b}R_{ijab}c(e^j\wedge e^a\wedge e^b)$ is zero. A cyclic change of the last three indices is even on the wedge, so renaming dummy indices and averaging gives
\begin{equation}
 3\mathcal T_i=\sum_{j,a,b}(R_{ijab}+R_{iabj}+R_{ibja})
                  c(e^j\wedge e^a\wedge e^b)=0.
\end{equation}
The last equality is the first Bianchi identity in its last-three-slot form, obtained from the usual first-three-slot identity by pair symmetry and antisymmetry. The remaining Clifford-linear part is
\begin{equation}
 \frac14\left(-\sum_{j,b}R_{ijjb}c_b
              +\sum_{j,a}R_{ijaj}c_a-\sum_{j,a}R_{ijaa}c_j\right).
\end{equation}
The third summand vanishes. Our Ricci convention gives
$\Ric_{ib}=\sum_jR_{jibj}=\sum_jR_{ijjb}$, whereas
$\sum_jR_{ijbj}=-\Ric_{ib}$. Thus the first two terms each contribute
$-\sum_b\Ric_{ib}c_b/4$, proving \eqref{eq:contractedspinexpanded}.
Contracting once more and using Ricci symmetry yields
\begin{equation}
 \sum_{i,j}c_ic_j\Omega^S_{ij}
 =-\frac12\sum_{i,b}\Ric_{ib}c_ic_b=\frac{s}{2}\id.
\end{equation}
Indeed $\sum\Ric_{ib}c_ic_b=\sum\Ric_{ib}(c_ic_b+c_bc_i)/2=-s\id$. Finally, moving $c_i$ through $C=c(Z)c(W)$ gives
$c_iC=Cc_i+2W_ic(Z)-2Z_ic(W)$. Right multiplication by $c_j\Omega^S_{ij}$ and summation now gives respectively $sC/2$, $-c(Z)c(\Ric^\sharp W)$ and $+c(W)c(\Ric^\sharp Z)$, exactly the last sandwich formula. Every sign follows from the specified four-slot curvature convention.

For use in the first evaluation, we record the scalar heat calculation in the body. If $L=D^2$, then $L=\nabla^{S*}\nabla^S+s/4$ and its first diagonal heat coefficient is $u_1=s/6-s/4=-s/12$. Thus for a zero-order insertion $A$ the coefficient of $t^{-m+1}$ in $\Tr(Ae^{-tL})$ is
$-(4\pi)^{-m}\int s\tr_S(A)\vol/12$.
Mellin transformation at the power $L^{-(m-1)}$ and the residue relation $\Res(P)=2\operatorname*{res}_{z=0}\Tr(PL^{-z})$ give
\begin{equation}
 \Res(AD^{-(n-2)})
 =-\frac{2(4\pi)^{-m}}{12\Gamma(m-1)}\int_Ms\tr_S(A)\vol
 =-\frac{\Cn(n-2)}{24}\int_Ms\tr_S(A)\vol.
 \label{eq:scalarinbody}
\end{equation}
Here $2(4\pi)^{-m}/\Gamma(m)=\Cn$ and $\Gamma(m)=(m-1)\Gamma(m-1)$. The factor two is the order of $L$: its degree $-n$ radial integral against $L^{-z}$ is $\int_1^\infty r^{-1-2z}dr=1/(2z)$. This calculation fixes the normalization used below.

\begin{theorem}[The corrected four-Clifford formulas]\label{thm:PQcorrect}
For $n=2m\ge4$, the first formula of \cite[Theorem 1.1]{LW} is
\begin{equation}
 P_D=-\frac{\Kn(n-2)}{24}\int_M
 s\bigl[g(X,Y)g(Z,W)-J(X,Y,Z,W)\bigr]\vol.
 \label{eq:Pcorrect}
\end{equation}
The second functional is instead
\begin{align}
 Q_D=\Kn\int_M\Big[&-\frac{n-2}{24}s\,g(X,Y)g(Z,W)
 +\frac{n-6}{24}s\,J(X,Y,Z,W)\nonumber\\
 &+\frac16(g\owedge\Ric)(X,Y,Z,W)\Big]\vol.
 \label{eq:Qcorrect}
\end{align}
Both expressions are integrated tensor densities; they contain no derivatives of the four argument fields.
\end{theorem}
\begin{proof}
For a zero-order insertion $A$, the heat coefficient calculation \eqref{eq:scalarinbody} yields
\[
 \Res(AD^{-(n-2)})=-\frac{\Cn(n-2)}{24}\int_M s\tr(A)\vol.
\]
Take $A=c(X)c(Y)c(Z)c(W)$ and use the four-generator trace to obtain \eqref{eq:Pcorrect}.

For $Q_D$ put $B=c(X)c(Y)$ and $C=c(Z)c(W)$. At a point choose a synchronous orthonormal frame. Expanding the differential operator $BDCD$, retaining the $D$ before $C$, gives
\[
 BDCD=B\sum_{i,j}c_iCc_j\nabla_i^S\nabla_j^S
       +B\sum_{i,j}c_i(\nabla_i^S C)c_j\nabla_j^S
\]
at that point, with frame-derivative corrections also of first order. All first-order terms have zero residue by Proposition \ref{prop:two}. The density of $Q_D$ is therefore
\begin{equation}
 \Cn\sum_{i,j}\tr\left(Bc_iCc_j
 \left(\frac{\Ric_{ij}}6+\frac{s\delta_{ij}}{24}
       +\frac{\Om^S_{ij}}2\right)\right).
 \label{eq:Qdensity}
\end{equation}
Combining the three identities in Lemma \ref{lem:PQclifford}, the endomorphism inside the trace after $B$ simplifies to
\begin{equation}
 \frac{6-n}{24}sC+\frac{s}{6}g(Z,W)\id
 +\frac16\bigl(c(W)c(\Ric^\sharp Z)-c(Z)c(\Ric^\sharp W)\bigr).
 \label{eq:Qmatrix}
\end{equation}
The scalar and four-Clifford traces are
\begin{align*}
 \tr B&=-r g(X,Y),\\
 \tr(BC)&=r\bigl[g(X,Y)g(Z,W)-J(X,Y,Z,W)\bigr],\\
 \tr\bigl(B(c(W)c(\Ric^\sharp Z)-c(Z)c(\Ric^\sharp W))\bigr)
 &=r(g\owedge\Ric)(X,Y,Z,W).
\end{align*}
In the last identity the two $g(X,Y)\Ric(Z,W)$ terms cancel, and the remaining four contractions are exactly \eqref{eq:KNproduct}. Substitution in \eqref{eq:Qmatrix} gives the coefficient $(6-n)/24-1/6=-(n-2)/24$ on $sg(X,Y)g(Z,W)$ and the other two coefficients in \eqref{eq:Qcorrect}.
\end{proof}

\begin{remark}[Checks and the missing Weyl information]
Setting all four inputs to the same unit field makes $J$ and $g\owedge\Ric$ vanish and agrees with the operator identity $Q_D(E,E,E,E)=\Res(D^{-(n-2)})$. Setting $Z=W=E$ in general agrees with $Dc(E)^2D=-D^2$. Clifford contractions in \eqref{eq:Qdensity} were also checked by exact rational arithmetic on all frame quadruples in dimensions four and six for algebraic curvature tensors; these checks supplement, rather than replace, the proof. Formulas \eqref{eq:Pcorrect}--\eqref{eq:Qcorrect} show that $P_D$, $Q_D$, and their permutations cannot distinguish algebraic curvature tensors with the same Ricci contraction. In particular the Weyl component is invisible to this family. The double commutator in Theorem \ref{thm:DiracR} removes that limitation.
\end{remark}

\subsection{A coefficient-by-coefficient audit of the corrected formula}\label{sec:Qannotations}
Here is an expanded ledger for \eqref{eq:Qmatrix}. Set $C=c(Z)c(W)$ and abbreviate $C_W=c(W)c(\Ric^\sharp Z)$, $C_Z=c(Z)c(\Ric^\sharp W)$. The three terms in \eqref{eq:Qdensity}, before the final trace with $B=c(X)c(Y)$, are respectively
\begin{align}
 \frac16\sum_{i,j}\Ric_{ij}c_iCc_j
 &=-\frac{s}{6}C-\frac13C_W+\frac13C_Z,\\
 \frac{s}{24}\sum_i c_iCc_i
 &=\frac{4-n}{24}sC+\frac{s}{6}g(Z,W)\id,\\
 \frac12\sum_{i,j}c_iCc_j\Om^S_{ij}
 &=\frac{s}{4}C+\frac12C_W-\frac12C_Z.
\end{align}
Thus the coefficient on $sC$ is $-1/6+(4-n)/24+1/4=(6-n)/24$; the coefficient on $C_W-C_Z$ is $-1/3+1/2=1/6$. After tracing with $B$, the scalar $g(X,Y)g(Z,W)$ coefficient becomes
\begin{equation}
 \frac{6-n}{24}-\frac16=-\frac{n-2}{24}.
\end{equation}
The $J$ coefficient is the negative of $(6-n)/24$, and the four Ricci cross-pairings have the same coefficient $1/6$. This equality of the cross-pairing coefficients is required by the algebraic curvature symmetries of $g\owedge\Ric$.

One may also verify the last trace directly, without notation for the product. Its first half is
\begin{align}
 r^{-1}\tr(B C_W)=g(X,Y)\Ric(Z,W)
 -g(X,W)\Ric(Y,Z)+g(Y,W)\Ric(X,Z),
\end{align}
whereas
\begin{align}
 r^{-1}\tr(B C_Z)=g(X,Y)\Ric(Z,W)
 -g(X,Z)\Ric(Y,W)+g(Y,Z)\Ric(X,W).
\end{align}
Subtracting cancels the terms proportional to $g(X,Y)$ and leaves exactly \eqref{eq:KNproduct}. A separate term pairing $g(X,Y)$ with $\Ric(Z,W)$ in the final formula would therefore fail this cancellation.

\begin{remark}[Differentiating a varying Clifford insertion]
In the operator expansion for $Q_D$, $D$ differentiates $c(Z)c(W)$. These derivative terms should not be dropped at the operator level. They are first-order insertions relative to the Bochner connection of $D^2$, so Proposition \ref{prop:two} shows that their residue is zero. The corrected formula is tensorial for that reason, even though the original differential expression does differentiate the input fields.
\end{remark}

\section{Spectral expressions for other curvature tensors}\label{sec:othercurv}
Here $n\ge4$. Write $P_D^-$ and $Q_D^-$ for half the difference obtained by exchanging $Z,W$. Multiplication by a scalar test function $a_0$ on the far left defines $P_{D,a_0}$ and $Q_{D,a_0}$; the calculations above are pointwise and apply without differentiating $a_0$.
\begin{corollary}[Ricci and Einstein curvature tensors]\label{cor:fourEinstein}
For every test function $a_0$,
\begin{align}
 P_{D,a_0}^-&=\frac{\Kn(n-2)}{24}\int_M a_0sJ\vol,\label{eq:Pminus}\\
 Q_{D,a_0}^- -\frac{n-6}{n-2}P_{D,a_0}^-
 &=\frac{\Kn}{6}\int_M a_0(g\owedge\Ric)\vol,\label{eq:Ricfour}\\
 Q_{D,a_0}^- -P_{D,a_0}^-
 &=\frac{\Kn}{6}\int_M a_0(g\owedge G)\vol,
 \qquad G=\Ric-\tfrac12sg.
 \label{eq:Einsteinfour}
\end{align}
\end{corollary}
\begin{proof}
The $g(X,Y)g(Z,W)$ terms are symmetric in $Z,W$; both $J$ and $g\owedge\Ric$ are antisymmetric. Antisymmetrizing Theorem \ref{thm:PQcorrect} gives \eqref{eq:Pminus} and \eqref{eq:Ricfour}. Subtracting $P^-$ rather than $(n-6)P^-/(n-2)$ adds $-\Kn sJ/6$. Since $g\owedge g=2J$, this is precisely \eqref{eq:Einsteinfour}. The tensor $g\owedge H$ has the algebraic curvature symmetries for every symmetric $H$; thus these antisymmetrized expressions, unlike the raw $P_D,Q_D$, define genuine algebraic curvature tensor densities.
\end{proof}

For compact notation denote by $\mathcal R_{a_0}$, $\mathcal L_{a_0}$ and $\mathcal S_{a_0}$ the localized integrated tensors recovered respectively by
\begin{equation}
 \mathcal R_{a_0}=-\frac{\mathcal U_D(a_0)}{4\Kn},\quad
 \mathcal L_{a_0}=\frac6{\Kn}\left(Q_{D,a_0}^--\frac{n-6}{n-2}P_{D,a_0}^-\right),\quad
 \mathcal S_{a_0}=\frac{24}{\Kn(n-2)}P_{D,a_0}^-.
 \label{eq:threebuilding}
\end{equation}
For the first expression use Corollary \ref{cor:oneforms} with the same four arguments as in $P,Q$; pair symmetry of $R$ aligns their order. Their respective densities are $R$, $g\owedge\Ric$, and $sJ$.

\begin{corollary}[Weyl, Schouten and concircular curvature]\label{cor:Weyl}
In our four-slot sign convention the Weyl tensor and the Schouten tensor are
\begin{align}
 \mathcal W&=R+\frac1{n-2}(g\owedge\Ric)
                 -\frac{s}{(n-1)(n-2)}J,\label{eq:Weyl}\\
 A&=\frac1{n-2}\left(\Ric-\frac{s}{2(n-1)}g\right).
\end{align}
Their spectral integrated realizations are
\begin{align}
 \mathcal W_{a_0}&=\mathcal R_{a_0}+\frac1{n-2}\mathcal L_{a_0}
                              -\frac1{(n-1)(n-2)}\mathcal S_{a_0},\\
 (g\owedge A)_{a_0}&=\frac1{n-2}\mathcal L_{a_0}
                              -\frac1{(n-1)(n-2)}\mathcal S_{a_0}.
\end{align}
The concircular and conharmonic tensors are recovered by
\begin{equation}
 \mathcal C_{a_0}=\mathcal R_{a_0}+\frac1{n(n-1)}\mathcal S_{a_0},
 \qquad
 \mathcal H_{a_0}=\mathcal R_{a_0}+\frac1{n-2}\mathcal L_{a_0}.
\end{equation}
\end{corollary}
\begin{proof}
For a symmetric $H$, contracting the first and fourth indices in \eqref{eq:KNproduct} gives $-(n-2)H-(\tr H)g$. The corresponding contraction of $J$ is $-(n-1)g$. These identities show that \eqref{eq:Weyl} is trace-free, and $\mathcal W=R+g\owedge A$ yields the Schouten decomposition in our convention. Substitution of the three recovered densities proves all spectral expressions. On constant sectional curvature $k$, our $R=-kJ$, $\Ric=(n-1)kg$, and $s=n(n-1)k$, so both $\mathcal W$ and $\mathcal C$ vanish as required. This also checks the signs.
\end{proof}

\begin{remark}[What the spectral realization means]
These are residue formulas for the tensors relative to the given Riemannian spectral triple. The metric $g$ and its Clifford action are part of this triple; no assertion that a single numerical residue reconstructs the metric is made. Theorem \ref{thm:DiracR} and the explicit correction to \cite{LW} are established here by calculation. These residue realizations are distinguished from connection-based constructions such as \cite{MRcurv}.
\end{remark}

\subsection{Expanded formulas for Weyl and the other curvature tensors}\label{sec:expandedWeyl}
To write everything as explicit residue combinations, let $U_{a_0}(X,Y,Z,W)$ denote the one-form version of \eqref{eq:Uoneforms} with $\alpha=X^\flat$, $\beta=Y^\flat$, $\gamma=Z^\flat$, $\delta=W^\flat$. Pair symmetry makes $-U_{a_0}/(4\Kn)$ the integrated $R(X,Y,Z,W)$. All $P^-,Q^-$ below have the same arguments and the left test function $a_0$.

The following identity is an entirely explicit Weyl residue formula:
\begin{equation}
 \boxed{\begin{aligned}
 \int_M a_0\mathcal W(X,Y,Z,W)\vol
 &=-\frac{U_{a_0}}{4\Kn}
   +\frac6{\Kn(n-2)}Q_{D,a_0}^-\\
 &\quad-\frac{6(n-5)}{\Kn(n-1)(n-2)}P_{D,a_0}^-.
 \end{aligned}}
 \label{eq:explicitWeylResidue}
\end{equation}
Indeed insert \eqref{eq:threebuilding} in \eqref{eq:Weyl}. The initially unsimplified coefficient on $P^-$ is
\begin{equation}
 -\frac{6(n-6)}{\Kn(n-2)^2}
 -\frac{24}{\Kn(n-1)(n-2)^2}
 =-\frac{6(n-5)}{\Kn(n-1)(n-2)},
\end{equation}
using $(n-6)(n-1)+4=(n-5)(n-2)$. Similarly, for the concircular tensor $\mathcal C=R+sJ/[n(n-1)]$ and the conharmonic tensor $\mathcal H=R+(g\owedge\Ric)/(n-2)$,
\begin{align}
 \int_M a_0\mathcal C\vol
 &=-\frac{U_{a_0}}{4\Kn}
   +\frac{24}{\Kn n(n-1)(n-2)}P_{D,a_0}^-,\\
 \int_M a_0\mathcal H\vol
 &=-\frac{U_{a_0}}{4\Kn}
   +\frac6{\Kn(n-2)}Q_{D,a_0}^-
   -\frac{6(n-6)}{\Kn(n-2)^2}P_{D,a_0}^-.
\end{align}
The arguments of the densities on the left are $(X,Y,Z,W)$. The Schouten curvature $g\owedge A$ is the Weyl formula without its $U$ term. The Einstein curvature is already the particularly simple difference \eqref{eq:Einsteinfour}.

To verify the trace removal in \eqref{eq:Weyl}, write $H$ for a symmetric tensor and compute its contraction directly:
\begin{align}
 \sum_a(g\owedge H)(e_a,Y,Z,e_a)
 &=H(Y,Z)+H(Y,Z)-nH(Y,Z)-(\tr_g H)g(Y,Z)\nonumber\\
 &=-(n-2)H(Y,Z)-(\tr_g H)g(Y,Z),\\
 \sum_aJ(e_a,Y,Z,e_a)&=-(n-1)g(Y,Z).
\end{align}
Thus the Ricci contraction of $R+(g\owedge\Ric)/(n-2)$ is $-sg/(n-2)$; subtracting $sJ/[(n-1)(n-2)]$ cancels it. This verifies the Weyl signs without appealing to a formula in a different slot convention. Antisymmetry and pair symmetry are immediate from \eqref{eq:KNproduct}; adding its three cyclic terms cancels them in pairs, proving the first Bianchi identity for $g\owedge H$.

In dimension four, the constants provide a useful short check:
\begin{equation}
 \Kn=\frac1{2\pi^2},\qquad
 \int_M a_0\mathcal W\vol=-\frac{U_{a_0}}{4\Kn}
                           +\frac3{\Kn}Q_{D,a_0}^-
                           +\frac1{\Kn}P_{D,a_0}^-.
\end{equation}
Substitution of $P^-_{a_0}=\Kn\int a_0sJ/12$ and $Q^-_{a_0}=\Kn\int a_0(-sJ/12+(g\owedge\Ric)/6)$ yields $R+(g\owedge\Ric)/2-sJ/6$, as expected.

\section{One-derivative insertions and skew torsion}
\subsection{The connection-defect principle}
\begin{proposition}\label{prop:defectM}
Let $L=\widetilde\nabla^*\widetilde\nabla+V$ and let $\nabla'=\widetilde\nabla+C$ be another connection on the same bundle. Then
\begin{equation}
 \Res(A\nabla'_ZL^{-m})=\Cn\int_M\tr(A C_Z)\,\vol.
 \label{eq:defectM}
\end{equation}
\end{proposition}
\begin{proof}
The insertion $A\widetilde\nabla_Z$ has zero residue by Proposition \ref{prop:two}. The remaining insertion $AC_Z$ is zero order; only $\sigma_{-n}(L^{-m})=|\xi|^{-n}\id$ contributes. Its sphere integral is $\sphere AC_Z$.
\end{proof}

\subsection{Torsion conventions and completion of the square}
Let $T$ be a real three-form and set
\begin{equation}
 \nabla^{T,TM}_ZW=\nabla^{LC}_ZW+\frac12T(Z,W,\cdot)^\sharp.
 \label{eq:torsionTM}
\end{equation}
Its torsion tensor is exactly $T(Z,W,\cdot)^\sharp$. The lifted connection and its Dirac operator are
\begin{align}
 \nabla^{T,S}_Z&=\nabla^S_Z+\frac14c(\iota_ZT),\\
 D_T&=D+B,\qquad B=\frac34c(T).
 \label{eq:DT}
\end{align}
Indeed the lift of the contorsion is $\frac18T_{Zab}c(e_a)c(e_b)=\frac14c(\iota_ZT)$, and $\sum_a c(e_a)c(\iota_{e_a}T)=3c(T)$. With $c(e_a)=i\gamma^a$, \eqref{eq:DT} becomes $B=-\frac i8 T_{abc}\gamma^a\gamma^b\gamma^c$, precisely the perturbation in \cite[equation (3.2)]{BDSZ}.

For a three-form, Clifford algebra gives
\begin{equation}
 \{c(e_a),c(T)\}=-2c(\iota_{e_a}T).
 \label{eq:antiT}
\end{equation}
At a normal frame, the first-order coefficient of $D_T^2$ relative to $D^2$ is $\{c(e_a),B\}\nabla^S_a$. Completing the connection square therefore gives the Bochner connection
\begin{equation}
 \widetilde\nabla_a=\nabla^S_a-\tfrac12\{c(e_a),B\}
 =\nabla^S_a+\tfrac34c(\iota_{e_a}T).
 \label{eq:bochnerT}
\end{equation}
This calculation uses only the first-order coefficient, so no curvature or derivative of $T$ can enter the connection defect. In particular,
\begin{equation}
 C_Z:=\nabla^{T,S}_Z-\widetilde\nabla_Z=-\tfrac12c(\iota_ZT).
 \label{eq:CT}
\end{equation}

\subsection{Detailed derivation of the lift and Bochner connection}\label{sec:torsionannotations}
The coefficient $3/4$ in \eqref{eq:DT}, and the discrepancy between the two connections, can be checked without using a torsion Lichnerowicz formula from another convention.

Let $K_ZW=\nabla^{T,TM}_ZW-\nabla^{LC}_ZW$. Then
\begin{equation}
 K_{Zab}=g(K_Ze_a,e_b)=\frac12T(Z,e_a,e_b).
\end{equation}
The spin representation of an infinitesimal skew matrix $K_{Zab}$ is $\frac14\sum_{a,b}K_{Zab}c_ac_b$. Thus the connection shift is
\begin{equation}
 \frac18\sum_{a,b}T(Z,e_a,e_b)c_ac_b=\frac14c(\iota_ZT).
\end{equation}
The $1/2$ in the component definition of a Clifford two-form is the reason that $1/8$ becomes $1/4$ here.

For every form $\alpha$ of degree $p$ the negative Clifford action obeys
\begin{align}
 c(v)c(\alpha)&=c(v^\flat\wedge\alpha)-c(\iota_v\alpha),\\
 c(\alpha)c(v)&=(-1)^p\bigl(c(v^\flat\wedge\alpha)+c(\iota_v\alpha)\bigr).
 \label{eq:wedgeContraction}
\end{align}
Both equations follow on a basis blade by separating whether $v$ belongs to the blade and then extend by linearity. For $p=3$ their sum is $-2c(\iota_vT)$, proving \eqref{eq:antiT}. Also
\begin{equation}
 \sum_a c_a c(\iota_{e_a}T)
 =c\left(\sum_a e_a^\flat\wedge\iota_{e_a}T\right)=3c(T).
 \label{eq:degreeOperator}
\end{equation}
The double contraction vanishes; the exterior degree operator in parentheses multiplies a three-form by three. Contracting the lifted connection with the $c_a$ proves $D_T=D+3c(T)/4$.

Now put $B=3c(T)/4$ and $U_i=c(\iota_{e_i}T)$. At a synchronous point, applying $(D+B)^2$ to a section gives
\begin{equation}
 D_T^2=D^2+\sum_i\{c_i,B\}\nabla_i^S
              +\sum_i c_i(\nabla_i^SB)+B^2.
 \label{eq:torsionSquareDetailed}
\end{equation}
The derivative of $B$ here is its induced endomorphism covariant derivative. Set $\widetilde\nabla_i=\nabla_i^S+\omega_i$. Expanding its connection Laplacian at this point gives
\begin{equation}
 \widetilde\nabla^*\widetilde\nabla
 =\nabla^{S*}\nabla^S-2\omega_i\nabla_i^S
                      -\nabla_i^S\omega_i-\omega_i\omega_i.
\end{equation}
Matching the first-order terms with \eqref{eq:torsionSquareDetailed} forces, rather than merely suggests,
\begin{equation}
 \omega_i=-\frac12\{c_i,B\}=\frac34U_i.
\end{equation}
The unique potential in $D_T^2=\widetilde\nabla^*\widetilde\nabla+V_T$ is correspondingly
\begin{equation}
 V_T=\frac s4\id+\sum_i c_i(\nabla_i^SB)+B^2
                   +\sum_i\nabla_i^S\omega_i+\sum_i\omega_i\omega_i.
 \label{eq:manifoldVT}
\end{equation}
This unsimplified formula is often safer than importing a simplified torsion potential with different lift conventions. It is an identity of endomorphisms; covariant divergences account for the synchronous-frame calculation globally. The one-derivative residue does not use $V_T$, because the coefficient $F$ in \eqref{eq:SIGMA} is zero. Derivatives of $T$, scalar curvature and quadratic torsion cannot therefore occur in \eqref{eq:Hfinal}.

Finally the inserted torsion connection has shift $U_i/4$, whereas the Bochner connection has shift $3U_i/4$. Their difference is $-U_i/2$. Using the degree-two trace in Section \ref{sec:traceannotations} gives the positive coefficient $\Kn/2$ of $T$ in \eqref{eq:Hfinal}. Using the Bochner connection instead sets this difference to zero. Using the Levi-Civita connection sets the difference to $-3U_i/4$ and gives the coefficient $3\Kn/4$. These three choices provide independent checks of the interpretation as a connection-defect functional.

\begin{remark}[The spectral torsion comparison, term by term]
The contraction \eqref{eq:degreeOperator} gives $\sum_i c_i\{c_i,B\}=-6B$. Thus the effective zero-order insertion in \eqref{eq:knownT}'s preceding calculation is $B-3B=-2B=-3c(T)/2$. Equation \eqref{eq:annSixTrace} then gives $-3\Kn T/2$. Dividing by the previously calculated $\Kn T/2$ gives exactly $-3$. The use of the degree-three, rather than degree-two, trace at this step is essential.
\end{remark}
\subsection{Evaluation of the one-derivative residue}

\begin{theorem}[Evaluation of the second functional]\label{thm:HT}
For \eqref{eq:torsionTM}--\eqref{eq:DT},
\begin{equation}
 \mathcal H_T(X,Y,Z)=\frac{\Kn}{2}\int_M T(X,Y,Z)\,\vol.
 \label{eq:Hfinal}
\end{equation}
The density is a three-form, and insertion of a scalar test function localizes it.
\end{theorem}
\begin{proof}
Apply \eqref{eq:defectM} with $L=D_T^2$ and $A=c(X)c(Y)$. Equations \eqref{eq:CT} and \eqref{eq:Ttrace} give
\begin{equation}
 \tr_S(A C_Z)=-\tfrac12\tr_S(c(X)c(Y)c(\iota_ZT))
 =\tfrac r2T(Z,X,Y)=\tfrac r2T(X,Y,Z).
\end{equation}
Multiplication by $\Cn$ proves \eqref{eq:Hfinal}.
\end{proof}

\subsection{Direct substitution in the SIGMA insertion formula}
The preceding proof can also be checked without naming the Bochner connection. At a Levi-Civita radial frame at the chosen point, put $U_a=c(\iota_{e_a}T)$ and $A=c(X)c(Y)$. For $D_T^2$, the first-order datum in \eqref{eq:SIGMA} is $S_a=\{c(e_a),B\}=-3U_a/2$. For $O=A\nabla_Z^{T,S}$,
\begin{equation}
 F^{ab}=0,\qquad G^a=A Z^a,\qquad H=\tfrac14A Z^aU_a.
\end{equation}
Thus \eqref{eq:SIGMA} reduces to
\begin{align}
 \Cn\tr(H+\tfrac12G^aS_a)
 &=\Cn\tr\bigl(A Z^a(\tfrac14U_a-\tfrac34U_a)\bigr)\\
 &=-\tfrac12\Cn\tr(Ac(\iota_ZT)).
\end{align}
This verifies both the sign and the factor in \eqref{eq:Hfinal}. The symbols $P,Q$ of $D_T^2$ are immaterial because $F=0$.

\subsection{Comparison with the known spectral torsion}
Define the functional of \cite{DSZtorsion}, in our Clifford convention, by
\begin{equation}
 \mathcal T_{D_T}(X,Y,Z)=\Res(c(X)c(Y)c(Z)D_TD_T^{-n}).
\end{equation}
For an arbitrary endomorphism $A_0$, the same first-order insertion calculation gives
\begin{equation}
 \Res(A_0D_TD_T^{-n})
 =\Cn\int_M\tr\left[A_0\left(B+\tfrac12\sum_a c(e_a)\{c(e_a),B\}\right)\right]\vol.
\end{equation}
For a three-form perturbation, the relevant contraction and trace are
\begin{equation}
 \sum_a c(e_a)\{c(e_a),B\}=-6B,\qquad
 \tr_S(c(X)c(Y)c(Z)c(T))=rT(X,Y,Z).
\end{equation}
Therefore
\begin{equation}
 \mathcal T_{D_T}(X,Y,Z)=-\frac{3\Kn}{2}\int_M T(X,Y,Z)\vol
 =-3\mathcal H_T(X,Y,Z).
 \label{eq:knownT}
\end{equation}
The coefficient in \eqref{eq:knownT} follows from the direct computation above and agrees with \cite[Proposition 3.11]{BDSZ} after the conversion $c=i\gamma$. We use that SIGMA coefficient consistently in numerical comparisons with the previously introduced functional. Hence the one-derivative torsion functional is a different presentation of the same three-form information. It therefore encodes the same three-form data.

\begin{remark}[Choice of connection and general torsion]
Using $\widetilde\nabla$ instead of $\nabla^{T,S}$ in \eqref{eq:Hintro} gives zero. Using $\nabla^S$ gives $3\Kn T/4$ as the integrated density. General metric torsion has vectorial, totally antisymmetric and Cartan components. The usual spin Dirac operator is blind to the Cartan component, and its self-adjointness imposes a condition on the vectorial component; see \cite{DSZtorsion}. Theorem \ref{thm:HT} assumes the three-form case and does not recover an arbitrary torsion tensor.
\end{remark}

\section{Ordered symbol calculus on a noncommutative torus}\label{sec:toruscalculus}
\subsection{Algebra, derivations and normalization}
Let $\theta$ be a real skew-symmetric $n\times n$ matrix. The smooth torus $\nc$ is generated by unitaries $U_j$ with
\begin{equation}
 U_jU_k=e^{2\pi i\theta_{jk}}U_kU_j.
\end{equation}
Its canonical trace is $\tau(\sum_k a_kU^k)=a_0$. Use the real-action derivations
\begin{equation}
 \dd_j(U^k)=ik_jU^k,\qquad \tau(\dd_j a)=0.
 \label{eq:realderiv}
\end{equation}
They are $*$-derivations and skew-adjoint on $L^2(\nc,\tau)$. Equivalently, the commonly used self-adjoint derivations are $\delta_j=-i\dd_j$. All coefficients act by left multiplication on the free right module $\nc^N$.

For matrix-valued classical symbols our composition convention is
\begin{equation}
 p\#q\sim\sum_\alpha\frac{1}{\alpha!i^{|\alpha|}}
 (\dd_\xi^\alpha p)(\dd^\alpha q),
 \label{eq:star}
\end{equation}
with the indicated order of multiplication. Define
\begin{equation}
 \Rtheta(P)=\int_{|\xi|=1}\tau\tr_N(p_{-n}(\xi))\,dS(\xi).
 \label{eq:NCres}
\end{equation}
This is the usual toroidal normalization. For $\theta=0$, our derivations correspond to the coordinate torus $\mathbb R^n/(2\pi\mathbb Z)^n$ and $\tau=(2\pi)^{-n}\int dx$. Thus \eqref{eq:NCres} agrees with \eqref{eq:resnormal}. Adding another $(2\pi)^{-n}$ to \eqref{eq:NCres} would incorrectly count this factor twice.

\subsection{Inverse-power symbols with noncommuting coefficients}
Consider
\begin{equation}
 L=-\sum_a\dd_a^2+b^a\dd_a+q,
 \qquad b^a,q\in M_N(\nc).
 \label{eq:flatL}
\end{equation}
The principal symbol is $|\xi|^2\id$. Since $m$ is an integer, $L^{-m}$ can be constructed by composing parametrices; a spectral cut is not required for this calculation.

\begin{lemma}\label{lem:NCinverse}
Write $\rho=|\xi|^2$. The first three homogeneous symbols of $L^{-m}$ are
\begin{align}
 r_{-n}&=\rho^{-m}\id,\\
 r_{-n-1}&=-mi b^a\xi_a\rho^{-m-1},\\
 r_{-n-2}&=-mq\rho^{-m-1}
 -\frac{m(m+1)}2b^ab^b\xi_a\xi_b\rho^{-m-2}\nonumber\\
 &\hspace{9mm}+m(m+1)(\dd_ab^b)\xi_a\xi_b\rho^{-m-2}.
 \label{eq:rinverse}
\end{align}
No commutation of $b^a,b^b$ is assumed.
\end{lemma}
\begin{proof}
The symbol of $L$ is $p_2=\rho$, $p_1=i b^a\xi_a$, $p_0=q$. Solve $p\#r=\id$ for a parametrix with $m=1$. At degrees $-2,-3,-4$ one obtains
\begin{align}
 r_{-2}&=\rho^{-1},\qquad r_{-3}=-i b^a\xi_a\rho^{-2},\\
 r_{-4}&=-q\rho^{-2}-b^ab^b\xi_a\xi_b\rho^{-3}
 +2(\dd_ab^b)\xi_a\xi_b\rho^{-3}.
\end{align}
For example the derivative contribution follows from
\begin{equation}
 \frac1i(\dd_{\xi_a}p_2)\dd_a r_{-3}
 =-2(\dd_ab^b)\xi_a\xi_b\rho^{-2},
\end{equation}
and multiplication by $-\rho^{-1}$ gives the last term in $r_{-4}$. For the induction, compose the formula for $L^{-m}$ with this inverse using \eqref{eq:star}. The coefficients of $q$, $b^ab^b$ and $\dd_ab^b$ become respectively
\begin{equation}
 -(m+1),\quad
 -\left(\frac{m(m+1)}2+1+m\right),\quad
 m(m+1)+2+2m.
\end{equation}
These are the coefficients in \eqref{eq:rinverse} with $m$ replaced by $m+1$. The $2m$ derivative contribution comes from differentiating the principal symbol of $L^{-m}$ in $\xi$ and the next symbol of $L^{-1}$ in the algebra derivations. Products keep their order throughout.
\end{proof}

\begin{theorem}[Ordered second-order insertion formula]\label{thm:NCmaster}
Let $L$ be \eqref{eq:flatL} and
\begin{equation}
 O=f^{ij}\dd_i\dd_j+g^i\dd_i+h,\qquad f^{ij}=f^{ji},
 \label{eq:Oflat}
\end{equation}
with arbitrary coefficients in $M_N(\nc)$. Then
\begin{align}
 \Rtheta(OL^{-m})=\sphere\,\tau\tr_N\biggl[&h+\tfrac12g^ib^i
 +f^{ii}\left(\tfrac12q-\tfrac14\dd_ab^a+\tfrac18b^ab^a\right)\nonumber\\
 &+\tfrac12 f^{ij}\dd_i b^j+\tfrac14 f^{ij}b^ib^j\biggr].
 \label{eq:NCmaster}
\end{align}
All repeated indices are summed, and the order displayed in \eqref{eq:NCmaster} is part of the formula.
\end{theorem}
\begin{proof}
The symbol of $O$ is $o_2=-f^{ij}\xi_i\xi_j$, $o_1=i g^i\xi_i$, $o_0=h$. The degree $-n$ symbol of $O\#L^{-m}$ is
\begin{equation}
 o_0r_{-n}+o_1r_{-n-1}+o_2r_{-n-2}
 +\frac1i(\dd_{\xi_a}o_2)\dd_a r_{-n-1}.
 \label{eq:fourcontributions}
\end{equation}
The other potential derivative terms vanish because $r_{-n}$ is constant in the algebra variables. Use the sphere moments
\begin{align}
 \int_{S^{n-1}}\xi_i\xi_j\,dS&=\frac{\sphere}{n}\delta_{ij},\\
 \int_{S^{n-1}}\xi_i\xi_j\xi_a\xi_b\,dS
 &=\frac{\sphere}{n(n+2)}
 (\delta_{ij}\delta_{ab}+\delta_{ia}\delta_{jb}+\delta_{ib}\delta_{ja}).
\end{align}
Since $n=2m$, $m/n=1/2$ and $m(m+1)/(n(n+2))=1/4$. The four summands in \eqref{eq:fourcontributions}, divided by $\sphere$ after integration, contribute
\begin{align}
 &h,\qquad \tfrac12g^ib^i,\\
 &\tfrac12 f^{ii}q+\tfrac18 f^{ii}b^ab^a+\tfrac14 f^{ij}b^ib^j
 -\tfrac14f^{ii}\dd_ab^a-\tfrac12f^{ij}\dd_i b^j,\\
 &f^{ij}\dd_i b^j,
\end{align}
respectively. To combine the two mixed quartic moments for $b^ib^j$, use symmetry of $f^{ij}$ and rename the indices; one does not commute $b^i$ with $b^j$. The sum is \eqref{eq:NCmaster}.
\end{proof}

\subsection{Full coefficient calculation for the ordered insertion}\label{sec:masterdetails}
We expand the sphere integrations in Theorem \ref{thm:NCmaster} before applying $\tau\tr_N$. Define a normalized sphere average
\begin{equation}
 \langle a\rangle=\sphere^{-1}\int_{|\xi|=1}a(\xi)\,dS.
\end{equation}
Thus $\langle h\rangle=h$, whereas the unnormalized integral of $h$ is $\sphere h$. This convention accounts for the common factor outside the final residue formula.

The left symbol is $o_2=-f^{ij}\xi_i\xi_j$, $o_1=ig^i\xi_i$, $o_0=h$. At degree $-n$ the complete list is
\begin{equation}
 o_0r_{-n}+o_1r_{-n-1}+o_2r_{-n-2}
 +i^{-1}(\partial_{\xi_a}o_1)\partial_ar_{-n}
 +i^{-1}(\partial_{\xi_a}o_2)\partial_ar_{-n-1}
 +\frac1{2i^2}(\partial_{\xi_a}\partial_{\xi_b}o_2)
                   \partial_a\partial_br_{-n}.
\end{equation}
Terms four and six vanish because $r_{-n}=\rho^{-m}\id$ has no algebra derivatives. Higher $\xi$ derivatives of $o_2$ vanish. This proves that the four surviving contributions are exhaustive.

On the sphere $\rho=1$. The first two averages are
\begin{equation}
 \langle o_0r_{-n}\rangle=h,\qquad
 \langle o_1r_{-n-1}\rangle
 =m g^ib^a\langle\xi_i\xi_a\rangle
 =\frac mn g^ib^i=\frac12g^ib^i.
\end{equation}
The $q$ contribution to the third is $mf^{ij}q\langle\xi_i\xi_j\rangle=f^{ii}q/2$. Its quadratic $b$ contribution is
\begin{align}
 &\frac{m(m+1)}2f^{ij}b^ab^b
       \langle\xi_i\xi_j\xi_a\xi_b\rangle\nonumber\\
 &\quad=\frac18\left(f^{ii}b^ab^a+f^{ij}b^ib^j+f^{ij}b^jb^i\right)
 =\frac18f^{ii}b^ab^a+\frac14f^{ij}b^ib^j.
\end{align}
For the last equality,
$\sum_{i,j}f^{ij}b^jb^i=\sum_{i,j}f^{ji}b^ib^j=\sum_{i,j}f^{ij}b^ib^j$.
Only dummy indices were renamed; the matrix products were never reversed. The derivative contribution to the third average is
\begin{align}
 &-m(m+1)f^{ij}(\partial_ab^b)
          \langle\xi_i\xi_j\xi_a\xi_b\rangle\nonumber\\
 &\quad=-\frac14\left(f^{ii}\partial_ab^a
                    +f^{ij}\partial_ib^j+f^{ij}\partial_jb^i\right)
 =-\frac14f^{ii}\partial_ab^a-\frac12f^{ij}\partial_ib^j.
\end{align}
The three delta pairings respectively set $(j=i,b=a)$, $(a=i,b=j)$ and $(b=i,a=j)$. Again $f^{ij}=f^{ji}$ combines the last two sums.

The remaining symbol derivative is essential. Symmetry of $f$ gives
\begin{equation}
 \partial_{\xi_a}o_2=-2f^{aj}\xi_j,\qquad
 \partial_ar_{-n-1}=-mi(\partial_ab^b)\xi_b\rho^{-m-1}.
\end{equation}
Their product with $1/i$ is
$2m f^{aj}(\partial_ab^b)\xi_j\xi_b\rho^{-m-1}$, whose sphere average is
$2mf^{aj}\partial_ab^j/n=f^{ij}\partial_ib^j$.
Combining it with the previous $-f^{ij}\partial_ib^j/2$ leaves the positive coefficient $+1/2$. The complete algebra-valued average is therefore the bracket in \eqref{eq:NCmaster}. Multiplication by $\sphere$ and application of $\tau\tr_N$ proves the theorem. There are no derivatives of $f,g,h$: the composition rule differentiates the left symbol in $\xi$ and the right symbol in the algebra directions.

\subsection{The inverse-power recursion in ordered form}\label{sec:inverseexpanded}
For clarity we record every nonzero source of the coefficients in Lemma \ref{lem:NCinverse}. Write $p_2=\rho$, $p_1=ib^a\xi_a$, $p_0=q$ and $s\sim s_{-2}+s_{-3}+s_{-4}$ for the inverse symbol. At degrees $0,-1,-2$ in $p\#s=1$,
\begin{align}
 \rho s_{-2}&=1,\nonumber\\
 \rho s_{-3}+p_1s_{-2}&=0,\nonumber\\
 \rho s_{-4}+p_1s_{-3}+qs_{-2}
       +i^{-1}(\partial_{\xi_a}p_2)\partial_as_{-3}&=0.
\end{align}
All omitted terms at these degrees contain an algebra derivative of $s_{-2}=\rho^{-1}$. Thus
\begin{align}
 s_{-3}&=-ib^a\xi_a\rho^{-2},\nonumber\\
 p_1s_{-3}&=b^ab^b\xi_a\xi_b\rho^{-2},\nonumber\\
 i^{-1}(\partial_{\xi_a}p_2)\partial_as_{-3}
 &=-2(\partial_ab^b)\xi_a\xi_b\rho^{-2}.
\end{align}
Multiplication of the last recursion by $-\rho^{-1}$ gives $s_{-4}$ in the lemma. The opposite signs of its quadratic and derivative terms are fixed by $i(-i)=1$ and $i^{-1}(-i)=-1$ respectively.

For induction denote the three orders of $L^{-m}$ by $R_0,R_1,R_2$. In $L^{-(m+1)}=L^{-m}L^{-1}$ the third order is exactly
\begin{equation}
 R_2s_{-2}+R_1s_{-3}+R_0s_{-4}
          +i^{-1}(\partial_{\xi_a}R_0)\partial_as_{-3}.
\end{equation}
The potential coefficients are $-m,-1$ from the first and third summands. The ordered $b^ab^b$ coefficients are $-m(m+1)/2,-m,-1$ from the first three. The derivative coefficients are $m(m+1),2,2m$ from the first, third and fourth, since
\begin{equation}
 i^{-1}(-2m\xi_a\rho^{-m-1})
        (-i\partial_ab^b\xi_b\rho^{-2})
 =2m(\partial_ab^b)\xi_a\xi_b\rho^{-m-3}.
\end{equation}
Their sums are $-(m+1)$, $-(m+1)(m+2)/2$ and $(m+1)(m+2)$, which are the next inductive values. The lower two orders are $\rho^{-m-1}$ and $-(m+1)ib^a\xi_a\rho^{-m-2}$. This establishes the complete recursion without a commutativity assumption. Ricci terms are absent because the principal metric is constant scalar, not because noncommutative connection curvature must vanish.

\section{Connection and Clifford formulas on the torus}
\subsection{Connection-compatible insertions}
Let
\begin{equation}
 \widetilde\nabla_i=\dd_i+\omega_i,\quad
 L=-\sum_i\widetilde\nabla_i^2+V,\quad
 \Om_{ij}=\dd_i\omega_j-\dd_j\omega_i+[\omega_i,\omega_j].
 \label{eq:NClaplace}
\end{equation}
The coefficients in \eqref{eq:flatL} are
\begin{equation}
 b^i=-2\omega_i,\qquad q=-\dd_i\omega_i-\omega_i\omega_i+V.
\end{equation}

\begin{corollary}\label{cor:NCconnection}
For any $A\in M_N(\nc)$ and fixed frame indices $i,j$,
\begin{align}
 \Rtheta(A\widetilde\nabla_iL^{-m})&=0,\label{eq:NCfirstzero}\\
 \Rtheta(A\widetilde\nabla_i\widetilde\nabla_jL^{-m})
 &=\frac{\sphere}{2}\tau\tr_N\bigl[A(\delta_{ij}V+\Om_{ij})\bigr].
 \label{eq:NCtwo}
\end{align}
\end{corollary}
\begin{proof}
For the first formula $g^a=A\delta_i^a$, $h=A\omega_i$ in \eqref{eq:NCmaster}; $h+g^ab^a/2=0$ exactly. For the second, take
\begin{align}
 f^{ab}&=\tfrac A2(\delta_i^a\delta_j^b+\delta_j^a\delta_i^b),\\
 g^a&=A(\omega_j\delta_i^a+\omega_i\delta_j^a),\qquad
 h=A(\dd_i\omega_j+\omega_i\omega_j).
\end{align}
The diagonal bracket in \eqref{eq:NCmaster} simplifies to $V/2$. The remaining terms simplify to
\begin{align}
 A\bigl(&\dd_i\omega_j+\omega_i\omega_j
 -\omega_j\omega_i-\omega_i\omega_j\nonumber\\
 &-\tfrac12(\dd_i\omega_j+\dd_j\omega_i)
 +\tfrac12(\omega_i\omega_j+\omega_j\omega_i)\bigr)
 =\tfrac12 A\Om_{ij}.
\end{align}
This proves \eqref{eq:NCtwo} without cyclic rearrangement of any coefficient.
\end{proof}

\begin{proposition}[Ordered connection defect]\label{prop:NCdefect}
Let $\nabla'_i=\widetilde\nabla_i+C_i$. Then
\begin{align}
 \Rtheta(A\nabla'_iL^{-m})&=\sphere\tau\tr_N(AC_i),\label{eq:NCdefectone}\\
 \Rtheta(A\nabla'_i\nabla'_jL^{-m})
 &=\sphere\tau\tr_N\left[A\left(\tfrac12\delta_{ij}V+\tfrac12\Om_{ij}
 +\dd_iC_j+[\omega_i,C_j]+C_iC_j\right)\right].
 \label{eq:NCdefecttwo}
\end{align}
\end{proposition}
\begin{proof}
The first identity follows from \eqref{eq:NCfirstzero} and the principal-symbol residue of $AC_iL^{-m}$. For the second, expand
\begin{align}
 \nabla'_i\nabla'_j
 =\widetilde\nabla_i\widetilde\nabla_j
 +C_j\widetilde\nabla_i+C_i\widetilde\nabla_j
 +\dd_iC_j+[\omega_i,C_j]+C_iC_j.
\end{align}
Both first-order summands have zero residue, even with their coefficient inserted to the left. The remaining zero-order terms contribute their principal symbols. This proves \eqref{eq:NCdefecttwo}.
\end{proof}

\subsection{Explicit cancellations and noncentral directions}\label{sec:connectiondetail}
The connection-compatible formulas can be verified entirely before taking traces. Expanding $(\partial_i+\omega_i)^2\psi$ gives
$\partial_i^2\psi+2\omega_i\partial_i\psi+(\partial_i\omega_i+\omega_i^2)\psi$, explaining the distinction between $q$ and $V$ in \eqref{eq:NClaplace}. The diagonal bracket in Theorem \ref{thm:NCmaster} is
\begin{equation}
 \frac q2-\frac14\partial_ab^a+\frac18b^ab^a
 =-\frac12\partial_a\omega_a-\frac12\omega_a^2+\frac V2
   +\frac12\partial_a\omega_a+\frac12\omega_a^2=\frac V2.
\end{equation}
For fixed $i,j$, the symmetric leading coefficient of
$A\widetilde\nabla_i\widetilde\nabla_j$ is
$f^{ab}=A(\delta_i^a\delta_j^b+\delta_j^a\delta_i^b)/2$, hence $\sum_af^{aa}=A\delta_{ij}$. Its remaining three contributions in the master formula are
\begin{align}
 \frac12g^ab^a&=-A\omega_j\omega_i-A\omega_i\omega_j,\nonumber\\
 \frac12f^{ab}\partial_ab^b&=-\frac A2(\partial_i\omega_j+\partial_j\omega_i),\nonumber\\
 \frac14f^{ab}b^ab^b&=\frac A2(\omega_i\omega_j+\omega_j\omega_i).
\end{align}
Adding $h=A(\partial_i\omega_j+\omega_i\omega_j)$ leaves
$A(\partial_i\omega_j-\partial_j\omega_i+\omega_i\omega_j-\omega_j\omega_i)/2=A\Omega_{ij}/2$.
This proves the compatible second-order formula by explicit ordered cancellation. Its symmetric half is $\sphere\delta_{ij}\tau\tr_N(AV)/2$, whereas its antisymmetric difference is $\sphere\tau\tr_N(A\Omega_{ij})$.

For an endomorphism $C_j$, the precise product rule is
\begin{equation}
 (\partial_i+\omega_i)(C_j\psi)
 =C_j(\partial_i+\omega_i)\psi
       +(\partial_iC_j+\omega_iC_j-C_j\omega_i)\psi.
\end{equation}
This accounts for the adjoint connection derivative in Proposition \ref{prop:NCdefect}. Subtracting its coefficient with $i,j$ interchanged yields
\begin{equation}
 \mathcal N_{ij}-\mathcal N_{ji}
 =\Omega_{ij}+\partial_iC_j-\partial_jC_i
       +[\omega_i,C_j]-[\omega_j,C_i]+[C_i,C_j],
\end{equation}
which is exactly the curvature of $\kappa_i=\omega_i+C_i$. This independent commutator check verifies the derivative and quadratic terms.

For directions $z^i,w^j$, the same product rule applies to their left multipliers. In particular
\begin{equation}
 \nabla'_i(w^j\psi)=w^j\nabla'_i\psi
                  +(\partial_iw^j+[\kappa_i,w^j])\psi.
\end{equation}
Multiplying by $z^i$ on the left produces \eqref{eq:arbitraryZW}. No derivative of $z^i$ occurs, since it is the outer left coefficient. The extra term vanishes if $C=0$ or if $\partial_iw^j+[\kappa_i,w^j]=0$. Centrality of $w^j$ alone removes only the commutator, not its derivative.

The correction is visible even at the commutative endpoint. Let $L=-\sum\partial_i^2$, $C_1=i$, $C_j=0$ for $j\ne1$, and use scalar rank one. Put $Z=\partial_1$, $W=w\partial_1$. Then
\begin{equation}
 \nabla'_Z\nabla'_W
 =w\partial_1^2+(2iw+\partial_1w)\partial_1-w+i\partial_1w.
\end{equation}
Only the zero-order coefficient contributes against $L^{-m}$. Thus the residue is $\sphere\tau(A(-w+i\partial_1w))$. If $w=e^{ix^1}$ and $A=w^*$, it is $-2\sphere$, whereas omitting the direction correction gives $-\sphere$. The issue is therefore an operator product rule coupled to connection mismatch, not solely a noncommutative phenomenon.

\subsection{What the flat connection example detects}\label{sec:bundleexampledetail}
For the flat scalar operator $D_0^2=-\sum\partial_i^2$, the inverse has symbol $\rho^{-m}$ modulo smoothing. A left-coefficient first- or second-order differential insertion has symbols of degrees $1-n$ or $2-n$, with no degree $-n$ term unless a zero-order coefficient is present. This is why \eqref{eq:flatvanish} remains zero for arbitrary algebra-valued $X,Y$. For $\partial_Z\partial_W=z^iw^j\partial_i\partial_j+z^i(\partial_iw^j)\partial_j$, both parts still have zero residue.

By contrast, a bundle connection commuting with the constant Clifford matrices has
\begin{equation}
 D_\omega^2=\sum_{i,j}c_ic_j\widetilde\nabla_i\widetilde\nabla_j
 =-\sum_i\widetilde\nabla_i^2+\sum_{i<j}c_ic_j\Omega_{ij}.
\end{equation}
For off-diagonal fixed indices and a Clifford-scalar test coefficient $h$,
\begin{equation}
 \Rtheta(h\widetilde\nabla_i\widetilde\nabla_jD_\omega^{-n})
 =\frac{r\sphere}{2}\tau\tr_k(h\Omega_{ij}),\qquad i\ne j.
\end{equation}
Although $\tau\tr_k\Omega_{ij}=0$ (derivatives and total commutators have zero trace), tests need not vanish. For instance take rank one, $\omega_1=0$, $\omega_2=i\varepsilon(u_1+u_1^*)$, other matrices zero. Then $\Omega_{12}=-\varepsilon(u_1-u_1^*)$ and
$\tau(\Omega_{12}^*\Omega_{12})=2\varepsilon^2$. The test $h=\Omega_{12}^*$ gives $r\sphere\varepsilon^2>0$. The base principal metric stays flat; this is auxiliary bundle curvature.

\subsection{Flat Dirac operator and a bundle connection}
Put $N=r$ and choose constant negative-square Clifford matrices $c_i$. The flat operator is
\begin{equation}
 D_0=\sum_i c_i\dd_i,\qquad D_0^2=-\sum_i\dd_i^2.
\end{equation}
For arbitrary algebra-valued frame coefficients $x^a,y^b$, define $c(X)=\sum_a x^ac_a$, $c(Y)=\sum_b y^bc_b$, with the indicated left order. Taking $A=c(X)c(Y)$ and $V=\omega=0$ in \eqref{eq:NCtwo} gives
\begin{equation}
 \Rtheta(c(X)c(Y)\dd_i\dd_jD_0^{-n})=0.
 \label{eq:flatvanish}
\end{equation}
The same is true for $\dd_Z\dd_W$, with $\dd_Z=z^i\dd_i$, $\dd_W=w^j\dd_j$: its extra term $z^i(\dd_iw^j)\dd_j$ has zero residue. Thus the naive flat counterpart of the first manifold functional contains no curvature.

To obtain nonzero connection curvature while retaining the flat principal metric, one can use a free module $\nc^r\otimes\mathbb C^k$ and $\omega_i\in M_k(\nc)$ commuting with all $c_a$. Then
\begin{equation}
 D_\omega=\sum_i c_i\widetilde\nabla_i,\qquad
 D_\omega^2=-\sum_i\widetilde\nabla_i^2+V_\omega,\qquad
 V_\omega=\sum_{i<j}c_ic_j\Om_{ij}.
\end{equation}
Equation \eqref{eq:NCtwo} gives an ordered trace pairing with $V_\omega$ and $\Om_{ij}$. These are the potential and bundle curvature of the selected connection. They are not the Riemann curvature of a curved base metric.

\subsection{Why arbitrary algebra-valued directions need extra terms}
Even if $\widetilde\nabla$ is the Bochner connection, coefficients $w^j\in\nc$ need not commute with $\omega_i$. Set $\nabla_Z=z^i\widetilde\nabla_i$, $\nabla_W=w^j\widetilde\nabla_j$. Direct expansion yields
\begin{align}
 \nabla_Z\nabla_W
 =z^iw^j\widetilde\nabla_i\widetilde\nabla_j
 +z^i\bigl(\dd_iw^j+[\omega_i,w^j]\bigr)\widetilde\nabla_j.
 \label{eq:directions}
\end{align}
Consequently the last term has zero residue by \eqref{eq:NCfirstzero}; the result is the pairing \eqref{eq:NCtwo} with $A z^iw^j$. For a different connection $\nabla'$, the same calculation gives
\begin{align}
 \Rtheta(A\nabla'_Z\nabla'_WL^{-m})
 =\sphere\tau\tr_N\bigl[&A z^iw^j\mathcal N_{ij}\nonumber\\
 &+A z^i(\dd_iw^j+[\kappa_i,w^j])C_j\bigr],
 \label{eq:arbitraryZW}
\end{align}
where $\kappa_i=\omega_i+C_i$ and
\begin{equation}
 \mathcal N_{ij}=\tfrac12\delta_{ij}V+\tfrac12\Om_{ij}
 +\dd_iC_j+[\omega_i,C_j]+C_iC_j.
\end{equation}
This extra term is generally nonzero. Frame-index formulas and ordered algebra-valued direction formulas must therefore be distinguished.

\subsection{Tested bundle curvature and the double-commutator limitation}\label{sec:NCcurvature}
The fixed-frame connection commutator has order zero even with noncommuting coefficients:
\begin{equation}
 [\widetilde\nabla_k,\widetilde\nabla_l]
 =\Om_{kl}=\dd_k\omega_l-\dd_l\omega_k+[\omega_k,\omega_l].
\end{equation}
Thus the leading inverse symbol alone gives a curvature pairing. This provides an ordered extension of the curvature insertion, rather than an unsupported extension of every commutative tensor identity.

\begin{proposition}[Fixed-frame noncommutative curvature pairing]\label{prop:NCcurv}
Let $L$ be as in \eqref{eq:NClaplace} on a free spinor module, and let $c_i$ be constant negative-square Clifford matrices. For $h\in\nc$,
\begin{equation}
 \Rtheta\bigl(h[c_a,c_b][\widetilde\nabla_k,\widetilde\nabla_l]L^{-m}\bigr)
 =\sphere\tau\tr_S\bigl(h[c_a,c_b]\Om_{kl}\bigr).
 \label{eq:NCcurvpair}
\end{equation}
If the connection curvature has spin degree two, written
\begin{equation}
 \Om_{kl}=\frac14\sum_{a,b}\mathcal R_{klab}c_ac_b,
 \qquad \mathcal R_{klab}=-\mathcal R_{klba}\in\nc,
 \label{eq:NCdegree2}
\end{equation}
then
\begin{equation}
 \Rtheta\bigl(h[c_a,c_b][\widetilde\nabla_k,\widetilde\nabla_l]L^{-m}\bigr)
 =-r\sphere\tau(h\mathcal R_{klab}).
 \label{eq:NCcurvrecover}
\end{equation}
The family of pairings for all $h$ determines each such coefficient.
\end{proposition}
\begin{proof}
The inserted commutator is multiplication by $\Om_{kl}$, so its product with $L^{-m}$ has degree $-n$ symbol $h[c_a,c_b]\Om_{kl}|\xi|^{-n}$. Integrating on the unit sphere proves \eqref{eq:NCcurvpair}. In \eqref{eq:NCdegree2} the coefficients commute with the constant Clifford matrices, though not necessarily with each other or with $h$. The Clifford trace in Lemma \ref{lem:trace} gives
\[
 \tr_S([c_a,c_b]\Om_{kl})=-r\mathcal R_{klab},
\]
without changing the order $h\mathcal R_{klab}$. If the pairing of a coefficient $t$ with every $h$ is zero, use $h=t^*$ and faithfulness of $\tau$ to conclude $t=0$.
\end{proof}

\begin{remark}[A necessary hypothesis on spin degree]
Formula \eqref{eq:NCcurvpair} holds for arbitrary matrix curvature. The representation \eqref{eq:NCdegree2} is an additional assumption, not an automatic consequence of degree-two connection coefficients over $\nc$. If coefficients of two bivectors do not commute, their matrix commutator can contain spin degrees zero and four as well. The tested spin trace always gives the degree-two projection; a general curvature endomorphism should be kept in the form \eqref{eq:NCcurvpair}. Even a degree-two curvature does not acquire pair symmetry, a Levi-Civita interpretation, or a first Bianchi identity merely from this residue formula. Those properties require additional metric and connection axioms.
\end{remark}

For clarity, the obstruction to copying Theorem \ref{thm:DiracR} can be seen already for the flat $D_0=\sum_i c_i\dd_i$. For $a\in\nc$,
\begin{equation}
 [D_0,a]=\sum_i c_i(\dd_i a),\qquad
 [D_0^2,a]=-2\sum_i(\dd_i a)\dd_i-\sum_i\dd_i^2a.
 \label{eq:NCsquarecomm}
\end{equation}
All coefficients act by left multiplication. The second-order part of the commutator of the last operators for $a,b$ is
\begin{equation}
 4\sum_{i,j}[\dd_i a,\dd_j b]\dd_i\dd_j,
 \label{eq:NCsecondobstruction}
\end{equation}
which need not vanish. Indeed the coefficient of $\dd_i^2$ is $4[\dd_i a,\dd_i b]$. Moreover,
\begin{equation}
 \tr_S([[D_0,a],[D_0,b]])=-r\sum_i[\dd_i a,\dd_i b],
 \label{eq:NCtraceobstruction}
\end{equation}
which need not be zero as an algebra element. Its canonical trace alone vanishes, but multiplying by another noncentral coefficient can produce a nonzero trace. Both cancellations in the classical proof can therefore fail. Equations \eqref{eq:NCsecondobstruction} and \eqref{eq:NCtraceobstruction} follow simply by expanding the two operator products and keeping coefficient order. The ordered insertion theorem above can evaluate the resulting terms in this flat-principal setting; they are not automatically a classical Riemann tensor. On a torus with scalar center, restricting all function arguments to the center would leave only constant functions and hence a trivial double-commutator functional.

\subsection{A noncentral curvature of pure Clifford degree two}\label{sec:spinexample}
Here is a concrete connection satisfying the degree-two assumption in Proposition \ref{prop:NCcurv}. On a noncommutative four-torus let $u_1u_2=q u_2u_1$ with $q\ne1$, put
\begin{equation}
 f=u_1+u_1^*,\qquad F=\partial_1f=i(u_1-u_1^*),\qquad
 \omega_2=f c_3c_4,\quad \omega_1=\omega_3=\omega_4=0.
\end{equation}
Both $f$ and $F$ are self-adjoint, while $\omega_2$ is skew-adjoint. At most one connection matrix is nonzero, so all cross-direction commutators of these matrices vanish. Differentiation gives
\begin{equation}
 \Omega_{12}=Fc_3c_4,\qquad \Omega_{21}=-Fc_3c_4,
 \qquad \Omega_{ij}=0\text{ for all other pairs}.
\end{equation}
No scalar or degree-four term occurs. In the normalization \eqref{eq:NCdegree2}, the independent coefficient is
$\mathcal R_{1234}=2F$, with $\mathcal R_{1243}=-2F$ and antisymmetry in the first pair. Indeed $\frac14(2Fc_3c_4-2Fc_4c_3)=Fc_3c_4$. The coefficient is noncentral, since
\begin{equation}
 [F,u_2]=i\bigl((1-q^{-1})u_1u_2-(1-q)u_1^*u_2\bigr)\ne0.
\end{equation}
The two Fourier monomials are distinct. However $\mathcal R_{3412}=0$, so pair-exchange symmetry fails. This is explicitly a bundle connection curvature, not an established Levi-Civita Riemann tensor.

For $L=-\sum\widetilde\nabla_i^2+\id$, Proposition \ref{prop:NCcurv} gives
\begin{equation}
 \Rtheta\bigl(h[c_3,c_4][\widetilde\nabla_1,\widetilde\nabla_2]L^{-2}\bigr)
 =-2r\nu_3\tau(hF).
\end{equation}
Since $\tau(F^2)=\tau(2-u_1^2-(u_1^*)^2)=2$, the choice $h=F$ gives $-4r\nu_3=-32\pi^2$. 

In this example $[\omega_2,c_a]$ is generally nonzero. Accordingly it uses the general connection-curvature formula, not the special square formula for $D_\omega$ in the preceding subsection, which assumes that the connection matrices commute with all Clifford generators.

\subsection{Ordered calculations: additional explanations}\label{sec:NCannotations}
The torus notation $\dd_i$ is a derivation of an algebra, not differentiation of an ordinary vector-valued coordinate function. Its operator Leibniz rule for a left coefficient $a$ is
\begin{equation}
 \dd_i\circ a=a\circ\dd_i+(\dd_i a).
 \label{eq:operatorLeibnizNC}
\end{equation}
In the last term parentheses indicate multiplication by the algebra element $\dd_i a$. This rule determines all lower-order coefficients in the ordered formulas.

\subsubsection{The first-order connection cancellation}
For $L=-\sum_i(\dd_i+\omega_i)^2+V$, expanding \eqref{eq:operatorLeibnizNC} gives
\[
 L=-\sum_i\dd_i^2-2\omega_i\dd_i
       -\sum_i\dd_i\omega_i-\sum_i\omega_i\omega_i+V.
\]
Thus $b^i=-2\omega_i$ and $q=-\dd_i\omega_i-\omega_i\omega_i+V$. An insertion $A\widetilde\nabla_j$ has $g^i=A\delta_j^i$ and $h=A\omega_j$. Its first-order residue is exactly
\begin{equation}
 \sphere\tau\tr\left(A\omega_j+\frac12A(-2\omega_j)\right)=0.
\end{equation}
The cancellation occurs as an algebra element before applying $\tau$, so a further noncentral coefficient $A$ does not spoil it. It uses the Bochner-compatible inserted connection, just as in the manifold calculation.

For a different connection $\nabla'_i=\widetilde\nabla_i+C_i$, expand
\begin{align}
 \nabla'_i\nabla'_j
 &=\widetilde\nabla_i\widetilde\nabla_j
   +C_j\widetilde\nabla_i+C_i\widetilde\nabla_j\nonumber\\
 &\quad+\dd_iC_j+[\omega_i,C_j]+C_iC_j.
 \label{eq:NCdefectExpandedNote}
\end{align}
The $[\omega_i,C_j]$ is required because $\widetilde\nabla_i\circ C_j$ must be put in left-coefficient order. The middle first-order terms have zero residue by the preceding cancellation, independently of $C_i,C_j$. This proves \eqref{eq:NCdefecttwo} directly from \eqref{eq:NCtwo}; the remaining zero-order products cannot be reordered.

If $Z=z^i e_i$ and $W=w^j e_j$, one further uses
\begin{equation}
 \nabla'_i\circ w^j=w^j\nabla'_i+\dd_iw^j+[\kappa_i,w^j],
 \qquad \kappa_i=\omega_i+C_i.
\end{equation}
Consequently
\begin{equation}
 \nabla'_Z\nabla'_W=z^iw^j\nabla'_i\nabla'_j
       +z^i(\dd_iw^j+[\kappa_i,w^j])\nabla'_j.
\end{equation}
In the compatible case $C_j=0$ the second term has zero residue; in the incompatible case its residue is
\begin{equation}
 \sphere\tau\tr\bigl(Az^i(\dd_iw^j+[\kappa_i,w^j])C_j\bigr).
\end{equation}
This is the additional term in \eqref{eq:arbitraryZW}. Noncentral direction coefficients cannot therefore be treated as ordinary tensor arguments by merely summing fixed-frame formulas.

\subsubsection{Why spin degree two is an extra hypothesis}
For algebra coefficients $u,v$ and constant Clifford matrices $P,Q$,
\begin{equation}
 [uP,vQ]=\frac12(uv+vu)[P,Q]+\frac12(uv-vu)\{P,Q\}.
\end{equation}
The anticommutator term is absent with commuting coefficients, but can create additional Clifford degrees in a noncommutative connection curvature. For instance,
\begin{equation}
 [u c_1c_2,v c_3c_4]=[u,v]c_1c_2c_3c_4,
 \qquad [u c_1c_2,v c_1c_2]=-[u,v]\id.
\end{equation}
Thus curvature can have degree four or zero even if the connection matrices are all degree two. The full spinor trace paired with $[c_a,c_b]$ selects the degree-two component, by \eqref{eq:gradeTrace}. It does not give unproved symmetries of the other components.

\subsubsection{A concrete failure of the classical double-commutator cancellation}
Take two generating unitaries with $UV=qVU$, $q\ne1$, and let $a=U$, $b=UV$. In their coordinate directions $\dd_1U=iU$, $\dd_1V=0$, and $\dd_2U=0$. Then
\begin{equation}
 [\dd_1a,\dd_1b]=-[U,UV]=-(1-q^{-1})U^2V\ne0.
\end{equation}
By \eqref{eq:NCsecondobstruction}, the coefficient of $\dd_1^2$ in $[[D_0^2,a],[D_0^2,b]]$ is $-4(1-q^{-1})U^2V$. It is genuinely second order, whereas the classical counterpart has order at most one. Equation \eqref{eq:NCtraceobstruction} also gives
\begin{equation}
 t:=\tr_S([[D_0,a],[D_0,b]])=r(1-q^{-1})U^2V\ne0.
\end{equation}
Although $\tau(t)=0$, taking a left test coefficient $h=t^*$ gives $\tau(ht)=\tau(t^*t)>0$. Thus trace cyclicity after applying $\tau$ cannot replace the vanishing of the spin trace as an algebra element. This example establishes the failure of the classical cancellation mechanism; it does not purport to define a noncommutative Levi-Civita Riemann tensor.

\subsubsection{Faithfulness and what the test pairings recover}
For the canonical Fourier trace, the monomials $U^k$ form an orthonormal system in $L^2(\nc,\tau)$. For $a=\sum_k a_kU^k$,
\begin{equation}
 \tau(a^*a)=\sum_k|a_k|^2.
\end{equation}
This identity proves faithfulness on the smooth algebra: if the left side is zero then every coefficient vanishes. For a matrix coefficient the same statement follows by summing the entrywise expressions. Hence the tests $\tau(ht)$ determine $t$ when all smooth $h$ are allowed; choosing $h=t^*$ already detects any nonzero difference of two candidate coefficients. Testing only $h=1$ determines its zero Fourier coefficient and can miss the entire nonconstant part. This distinction is used both in \eqref{eq:NCcurvrecover} and \eqref{eq:NClocal}.

\begin{remark}[Normalization at the commutative endpoint]
With $\dd_j(U^k)=ik_jU^k$, the coordinate torus has period $2\pi$ in each direction. At $\theta=0$, $\tau(a)=(2\pi)^{-n}\int_{[0,2\pi]^n}a(x)\,dx$. Therefore $\sphere\tau\tr$ already equals $\Cn\int\tr$. A further $(2\pi)^{-n}$ would be an error. In dimension four $r\sphere/2=4\pi^2$ in a $\tau$ pairing, while the corresponding normalized manifold coefficient $\Kn/2$ is $1/(4\pi^2)$ in an ordinary volume integral; the difference is exactly the normalized torus volume factor.
\end{remark}

\section{Torsion perturbations on the noncommutative torus}
\subsection{Algebra-valued three-form coefficients}
Let $T_{abc}\in\nc$ be smooth, self-adjoint and totally antisymmetric in their indices. They commute with the constant Clifford matrices but need not commute with one another. Define
\begin{align}
 t=c(T)&=\tfrac16T_{abc}c_ac_bc_c,\qquad U_i=c(\iota_iT)=\tfrac12T_{iab}c_ac_b,\\
 \kappa_i&=\tfrac14U_i,\qquad \nabla^T_i=\dd_i+\kappa_i,\\
 D_T&=\sum_i c_i\nabla^T_i=D_0+\tfrac34t.
 \label{eq:NCT}
\end{align}
These define a spectral/connection torsion model on the free spinor module. No uniqueness theorem for a bimodule Levi-Civita connection is assumed. The operator $D_T$ is self-adjoint on the flat Dirac domain: the bounded perturbation $t$ is self-adjoint because a product of three distinct negative-square skew-adjoint Clifford matrices is self-adjoint.

Expanding its square in the order of \eqref{eq:flatL} gives
\begin{equation}
 b^i=\{c_i,\tfrac34t\}=-\tfrac32U_i,
 \qquad q=\tfrac34c_a\dd_at+\tfrac9{16}t^2.
\end{equation}
Hence its Bochner decomposition is
\begin{align}
 \omega_i&=\tfrac34U_i,\qquad C_i=\kappa_i-\omega_i=-\tfrac12U_i,\label{eq:NCbochner}\\
 V_T&=\tfrac34c_a\dd_at+\tfrac9{16}t^2
 +\tfrac34\dd_aU_a+\tfrac9{16}U_aU_a.
 \label{eq:NCV}
\end{align}
These identities involve no commutation of the $T$ coefficients.

\begin{theorem}[Ordered torus torsion recovery]\label{thm:NCtorsion}
Let $X=(x^a)$, $Y=(y^b)$, $Z=(z^i)$ have arbitrary smooth coefficients in $\nc$, and set $\nabla^T_Z=z^i\nabla^T_i$. Then
\begin{equation}
 \Rtheta\bigl(c(X)c(Y)\nabla^T_ZD_T^{-n}\bigr)
 =\frac{r\sphere}{2}\sum_{a,b,i}\tau(x^ay^bz^iT_{iab}).
 \label{eq:NCtorsion}
\end{equation}
For fixed frame indices $a,b,i$ and any test coefficient $f\in\nc$,
\begin{equation}
 \Rtheta(f c_ac_b\nabla^T_iD_T^{-n})
 =\frac{r\sphere}{2}\tau(fT_{iab}).
 \label{eq:NClocal}
\end{equation}
All the coefficients $T_{iab}$ are determined by the collection of these test pairings.
\end{theorem}
\begin{proof}
Use \eqref{eq:NCdefectone} with $A=c(X)c(Y)z^i$ and $C_i=-U_i/2$. The spin trace identity is an identity over the coefficient algebra,
\begin{equation}
 \tr_S(c_ac_bU_i)=-rT_{iab}.
\end{equation}
The product preceding $T_{iab}$ remains $x^ay^bz^i$; no cyclic permutation is needed. This proves \eqref{eq:NCtorsion} and \eqref{eq:NClocal}. If the test pairing of a coefficient $T_{iab}$ vanishes for every $f$, choose $f=T_{iab}^*$. Faithfulness of $\tau$ implies $T_{iab}=0$, proving separation by test coefficients.
\end{proof}

\subsection{Failure of naive alternation}
Formula \eqref{eq:NCtorsion} is alternating in its frame indices, but is generally not alternating under exchange of arbitrary algebra-valued arguments. For example, take $X=a e_1$, $Y=b e_2$, $Z=e_3$, and $T_{312}=t$. Then
\begin{equation}
 \mathcal H_{\theta,T}(X,Y,Z)+\mathcal H_{\theta,T}(Y,X,Z)
 =\frac{r\sphere}{2}\tau([a,b]t).
 \label{eq:alternationobstruction}
\end{equation}
For noncommuting self-adjoint $a,b$, take the self-adjoint $t=i[a,b]$. Since $[a,b]^*=-[a,b]$,
\begin{equation}
 \tau([a,b]t)=i\tau([a,b]^2)=-i\tau([a,b]^*[a,b])\ne0.
\end{equation}
Thus cyclicity of $\tau$ does not justify interchanging adjacent factors. One can instead restrict argument coefficients to the center, or explicitly antisymmetrize the entire functional. Either change must be stated as part of the definition.

\begin{proposition}[Centrality from restricted alternation]\label{prop:centrality}
Fix distinct frame indices $a,b,i$ and put $t=T_{iab}$. The restriction of \eqref{eq:NCtorsion} to $X=u e_a$, $Y=v e_b$, $Z=e_i$ is alternating in $X,Y$ for all $u,v\in\nc$ if and only if $t$ belongs to the center of $\nc$.
\end{proposition}
\begin{proof}
The condition is $\tau([u,v]t)=0$ for every $u,v$. By cyclicity this is equivalent to $\tau(u[v,t])=0$ for every $u,v$. Fix $v$ and take $u=[v,t]^*$. Faithfulness yields $[v,t]=0$. Conversely, centrality gives $\tau([u,v]t)=\tau([ut,v])=0$.
\end{proof}

\begin{theorem}[Obstruction to universal ordinary alternation]\label{thm:universalalt}
Assume $\nc$ is noncommutative and has center $\mathbb C\id$. For $n\geq4$, the functional \eqref{eq:NCtorsion} is alternating under exchange of its first two arguments for every algebra-valued $X,Y,Z$ if and only if $T=0$.
\end{theorem}
\begin{proof}
Suppose this universal alternation holds. Taking $Z=e_i$ and applying Proposition \ref{prop:centrality} shows that every independent component $T_{iab}$ is a scalar $\lambda\id$. If some component is nonzero, choose $X=u e_a$, $Y=v e_b$, $Z=w e_i$ with $[u,v]\ne0$. The sum under exchange of $X,Y$, divided by the nonzero factor $r\sphere/2$, equals
\begin{equation}
 \lambda\tau([u,v]w).
\end{equation}
Choose $w=[u,v]^*$. Then $\tau([u,v][u,v]^*)>0$, contradicting alternation. All components must vanish. The converse is immediate.
\end{proof}

\begin{remark}
Theorem \ref{thm:universalalt} concerns ordinary alternation of this particular ordered numerical functional with unrestricted left coefficient arguments. It is not an obstruction to defining torsion in a bimodule differential calculus, or to braided notions of alternation. Even a nonzero constant central torsion component fails the stated universal requirement on a noncommutative torus. Its distinction from tensoriality obstructions for spectral Einstein functionals \cite{Impediment} is therefore essential.
\end{remark}

\subsection{The torsion-deformed two-derivative formula}
As an additional extension of the first functional, use the same $D_T$ in both the connection insertion and the inverse power. For fixed frame directions,
\begin{equation}
 \Rtheta(A\nabla_i^T\nabla_j^TD_T^{-n})
 =\sphere\tau\tr_S(A\mathcal N^T_{ij}),
 \label{eq:NCTwoT}
\end{equation}
where the completely ordered coefficient is
\begin{equation}
 \boxed{\mathcal N^T_{ij}
 =\tfrac12\delta_{ij}V_T-\tfrac18\dd_iU_j-\tfrac38\dd_jU_i
 +\tfrac5{32}U_iU_j+\tfrac3{32}U_jU_i.}
 \label{eq:NCTexplicit}
\end{equation}
Indeed \eqref{eq:NCdefecttwo} and \eqref{eq:NCbochner} give
\begin{align}
 \tfrac12\Om_{ij}&=\tfrac38(\dd_iU_j-\dd_jU_i)+\tfrac9{32}[U_i,U_j],\\
 \dd_iC_j+[\omega_i,C_j]&=-\tfrac12\dd_iU_j-\tfrac38[U_i,U_j],\\
 C_iC_j&=\tfrac14U_iU_j.
\end{align}
Combining them yields \eqref{eq:NCTexplicit}. Taking $A=c(X)c(Y)$ gives the desired Clifford-weighted extension; leaving the Clifford products in this form preserves every ordered quadratic torsion term explicitly. No implicit ``commutative contraction'' is used.

Its antisymmetric component in the derivative frame indices is particularly simple:
\begin{equation}
 \mathcal N^T_{ij}-\mathcal N^T_{ji}
 =\tfrac14(\dd_iU_j-\dd_jU_i)+\tfrac1{16}[U_i,U_j]
 =\Om^\kappa_{ij}.
 \label{eq:NCantisymT}
\end{equation}
This is the curvature of the inserted connection $\nabla^T$, not of the Bochner connection. It also follows directly by taking the operator commutator and observing that $\Om^\kappa_{ij}D_T^{-n}$ contributes only its leading symbol. For arbitrary algebra-valued $Z,W$, one must additionally include the second line of \eqref{eq:arbitraryZW}.

\subsection{Ordered quadratic terms and a commutator check}\label{sec:torsionquadraticdetail}
For the torsion model, the general defect coefficient is
\begin{equation}
 \mathcal N^T_{ij}=\frac12\delta_{ij}V_T+\frac12\Omega_{ij}
             +\partial_iC_j+[\omega_i,C_j]+C_iC_j.
\end{equation}
The three non-potential pieces give respectively the derivative coefficients
$(3/8,-3/8)$ and $(-1/2,0)$ on $(\partial_iU_j,\partial_jU_i)$, so their sum is $(-1/8,-3/8)$. The commutator coefficient is $9/32-3/8=-3/32$. Combining it with $C_iC_j=U_iU_j/4$ gives
\begin{equation}
 -\frac3{32}(U_iU_j-U_jU_i)+\frac8{32}U_iU_j
 =\frac5{32}U_iU_j+\frac3{32}U_jU_i.
\end{equation}
These two products differ even for commuting torsion coefficients, because of their Clifford factors. For $T_{123}=s$, one has
$U_1U_2=-s^2c_1c_2$ but $U_2U_1=s^2c_1c_2$.

The symmetric half of the coefficient is
\begin{equation}
 \frac{\mathcal N^T_{ij}+\mathcal N^T_{ji}}2
 =\frac12\delta_{ij}V_T-\frac14(\partial_iU_j+\partial_jU_i)
                  +\frac18(U_iU_j+U_jU_i).
\end{equation}
The antisymmetric difference is \eqref{eq:NCantisymT}; its half, rather than its full difference, is the conventionally normalized antisymmetric part. The curvature coefficient in the full difference is one, as can also be checked without expansion:
\begin{equation}
 \Rtheta(A[\nabla_i^T,\nabla_j^T]D_T^{-n})
 =\sphere\tau\tr_S(A\Omega^\kappa_{ij}).
\end{equation}
For arbitrary direction coefficients the complete specialization of \eqref{eq:arbitraryZW} is
\begin{equation}
 \begin{split}
 \Rtheta(A\nabla_Z^T\nabla_W^TD_T^{-n})
 =\sphere\tau\tr_S\biggl[&Az^iw^j\mathcal N^T_{ij}\\
 &-\frac12Az^i\left(\partial_iw^j+\frac14[U_i,w^j]\right)U_j\biggr].
 \end{split}
\end{equation}
This displays the direction correction instead of absorbing it into a putative tensor coefficient.

\subsection{Testing a four-dimensional noncentral component}\label{sec:torsionexampledetail}
For the example in \eqref{eq:example4}, the six nonzero coefficients are
$T_{123}=T_{231}=T_{312}=t$ and $T_{132}=T_{213}=T_{321}=-t$.
In the unrestricted sum for $c(T)$ all six products equal $tc_1c_2c_3$, because the coefficient and Clifford permutation signs cancel. Likewise
\begin{equation}
 U_3=\frac12(T_{312}c_1c_2+T_{321}c_2c_1)=tc_1c_2.
\end{equation}
The inserted-to-Bochner connection difference is $C_3=-tc_1c_2/2$. Therefore
\begin{equation}
 fc_1c_2C_3=-\frac12ft(c_1c_2)^2=\frac12ft\id.
\end{equation}
The four-dimensional spin trace is $2ft$, and $\nu_3=2\pi^2$, giving exactly $4\pi^2\tau(ft)$ in \eqref{eq:example4}.

If $t=u_1+u_1^*$, then $\tau(t)=0$ but
\begin{equation}
 t^2=u_1^2+2+(u_1^*)^2,\qquad \tau(t^2)=2.
\end{equation}
Thus the tests $f=1$ and $f=t$ give respectively zero and $8\pi^2$. In general $\partial_1t=i(u_1-u_1^*)\ne0$; no constant-torsion assumption has been made. Derivatives of torsion do not enter this one-derivative residue because the relevant inverse order sees only the connection defect. When $u_1u_2=q u_2u_1$ with $q\ne1$,
$[t,u_2]=(1-q^{-1})u_1u_2+(1-q)u_1^*u_2\ne0$, so this is a noncentral component. A single vanishing unweighted residue measures only zero mean, not zero torsion.

For completeness the general first-pair alternation defect follows directly from ordered coefficient sums:
\begin{equation}
 \mathcal H_{\theta,T}(X,Y,Z)+\mathcal H_{\theta,T}(Y,X,Z)
 =\frac{r\sphere}{2}\sum_{a,b,i}\tau([x^a,y^b]z^iT_{iab}).
\end{equation}
To obtain it rename $a,b$ in the second functional and use $T_{iba}=-T_{iab}$. This renaming changes $y^ax^b$ to $y^bx^a$; it does not change their order. The identity
$\tau([u,v]t)=\tau(u[v,t])$ follows from $\tau(vut)=\tau(utv)$, a cyclic move, not an adjacent transposition. Taking $u=[v,t]^*$ in a vanishing pairing yields $\tau([v,t]^*[v,t])=0$ and hence centrality by faithfulness, as used in Proposition \ref{prop:centrality}. Allowing a third weight $w$ instead gives $\tau([u,v]wt)$; even a nonzero scalar $t$ cannot annihilate every such test in a noncommutative algebra. This explains the distinction between restricted and universal alternation in Theorem \ref{thm:universalalt}.

\subsection{A four-dimensional noncentral example}
Take $n=4$, and let the only independent torsion component be $T_{123}=t=t^*\in\nc$. Then
\begin{equation}
 c(T)=t c_1c_2c_3,\quad
 U_1=t c_2c_3,\quad U_2=-t c_1c_3,\quad U_3=t c_1c_2,\quad U_4=0.
\end{equation}
Since $r=4$ and $\nu_3=2\pi^2$,
\begin{equation}
 \Rtheta(f c_1c_2\nabla^T_3D_T^{-4})=4\pi^2\tau(ft).
 \label{eq:example4}
\end{equation}
Thus $f=t$ produces $4\pi^2\tau(t^2)>0$ for any nonzero self-adjoint $t$. Let $u_1$ denote the first generating unitary of the torus algebra (as distinct from the contraction endomorphisms $U_i$ above). A noncentral example, when this generator is noncentral, is $t=u_1+u_1^*$; in the canonical Fourier normalization $\tau(t^2)=2$. A constant test coefficient $f=1$ gives zero for this example, whereas the localized test detects it. This illustrates why the family of test pairings, rather than the single unlocalized residue, is needed for recovery.

\section{Conclusions and scope}
The classical Dirac double commutator yields a localized full-curvature pairing, and its extension through scalar one-form presentation errors is well defined after taking the residue. Four-Clifford insertions supply complementary scalar and Ricci data. Their corrected coefficients, together with the full-curvature pairing, produce the explicit Weyl, Schouten, Einstein, concircular and conharmonic expressions proved above. The torsion insertion detects the inserted-to-Bochner connection difference and realizes the established three-form information through a covariant derivative.

For general spectral triples the exact last slots of Definition \ref{def:partial} make the partial pairing intrinsic without assuming presentation descent. The flat torus counterexamples show why this is separate from being a module tensor or the curvature of a specified Levi-Civita connection. The fixed-frame connection-curvature pairing identifies its tested Clifford-degree-two coefficients; further geometric axioms are needed to identify these with a Levi-Civita curvature tensor. 

The ordered torus formulas require a free module, smooth left coefficients, flat scalar principal symbol, and the residue normalization \eqref{eq:NCres}. Arbitrary coefficient directions require the correction \eqref{eq:arbitraryZW}. The ordinary alternation obstruction concerns this particular ordered numerical torsion functional, not braided or bimodule definitions of torsion. Boundary problems, projective modules and noncentral curved principal metrics are outside the calculations proved here.

All results concern classical or toroidal pseudodifferential residue calculi, or the explicitly assumed abstract residue domain in Definition \ref{def:partial}. They do not assert that an arbitrary finitely summable triple automatically supports every displayed composition.

\section*{Acknowledgements}
The author was supported by Science and Technology Development Plan Project of Jilin Province, China: No.20260102245JC and NSFC No.11771070.

\end{document}